\documentclass[11pt]{article}
\usepackage{amsmath,amsthm}
\usepackage{pstricks,pst-node,pst-tree}
\usepackage{latexsym}
\usepackage{amssymb,amscd}
\usepackage{graphics}   
\usepackage{url}
\numberwithin{equation}{section}
\let \:=\colon

\begin{document}
\title{Rational curves of degree $10$ on a general quintic threefold \thanks{$2000$ {\it Mathematics Subject Classification}.
Primary 14J30; Secondary 13P10, 14H45.}}
\author{Ethan Cotterill \thanks{Partially supported by an NSF graduate fellowship.}}
\vspace{-.5cm}
\date{December, 2004}
\thispagestyle{empty}
\maketitle
\begin{abstract}
We prove the ``strong form'' of the Clemens conjecture in degree
$10$. Namely, on a general quintic threefold $F$ in $\mathbb{P}^4$,
there are only finitely many smooth rational curves of degree $10$, and
each curve $C$ is embedded in $F$ with normal bundle $\mathcal{O}(-1) \oplus
\mathcal{O}(-1)$. Moreover, in degree $10$, there are no singular,
reduced, and irreducible rational
curves, nor any reduced, reducible, and connected curves with rational
components on $F$.
\end{abstract}

\section* {Introduction}
Almost twenty years ago, Clemens conjectured that the number of smooth
rational curves of fixed degree on a general quintic threefold in
$\mathbb{P}^4$ is finite. In 1986, Sheldon Katz \cite{Katz} proved Clemens' conjecture
in degrees at most $7$. His method of proof was in two
steps. First he used a
deformation-theoretic argument to show that

{\it Clemens' conjecture holds for smooth rational curves
  $C$  of degree $d$ provided the incidence scheme 
\vspace{-.25cm}
\begin{equation}
\Phi_d:=\{
(C,F)\mid\text{$F$ a quintic containing $C$}\}
\notag 
\vspace{-.25cm}
\end{equation}
is irreducible.}

He then proved that $\Phi_d$ is irreducible for all $d \leq
7$, by arguing that the
fibres of its projection onto the Hilbert scheme of smooth degree-$d$
curves are equidimensional projective spaces. Indeed, equidimensionality follows immediately from a
vanishing result for ideal-sheaf cohomology of Gruson, Lazarsfeld, and
Peskine.

Building on Katz's work, Johnsen and Kleiman \cite{Klei} next showed
that the only reduced connected curves of degree $d$ at most 9 on $F$
with rational components are irreducible and either smooth or six-nodal
plane quintics. (In the meantime, Vainsencher \cite{Va} had
discovered that there are 17,601,000
six-nodal plane quintics on $F$. For a proof of the latter
assertion, see \cite{Piene}.) To establish that the generic
threefold $F$ contains curves of the asserted sort, Johnsen and Kleiman first observe
that the space 
\vspace{-.2cm}
\begin{equation*}
M_d:=\{\text{degree-$d$ morphisms }
f\:\mathbb{P}^1 \rightarrow \mathbb{P}^4\}
\vspace{-.2cm}
\end{equation*}
is stratified by locally-closed subsets
$M_{d,i}$, where $i:=h^1(\mathcal{I}_{C/\mathbb{P}^4}(5))$ and $C$ is
the image of $f$.
Pulling back by $\pi_d:\Phi_d \rightarrow M_d$, they obtain a
corresponding stratification of the incidence scheme $\Phi_d$ into loci
$I_{d,i}$, and they show that the projection $I_{d,i}
\rightarrow M_{d,i}$ is dominant exactly over $M_{d,0}$, where all the
fibres are equidimensional. 

Johnsen and Kleiman apply the same result of
Gruson--Lazarsfeld--Pes\-kine used by Katz, as well as a related result
of d'Almeida, to show that whenever the fibre dimension
$h^0(\mathcal{I}_{C/\mathbb{P}^4}(5))-1$ is not the expected one, then the curve $C$
necessarily admits a highly-incident secant line. Smooth curves
with high-incidence secants comprise a locally closed locus of high
codimension inside the Hilbert scheme; Johnsen and Kleiman settle Clemens' conjecture
for smooth curves by obtaining appropriate bounds on the dimensions of
the fibres $\pi_d^{-1}(C)$ over this locus.

To handle singular irreducible and reduced rational curves $C$, Johnsen and Kleiman use the
same basic approach, but are forced to work
harder, because for singular curves the arithmetic genus $g(C)$
varies over a finite nonzero range. To bound this
variation, they apply a classical result of
Castelnuovo--Halphen giving the maximal possible arithmetic genus of each curve in
terms of the dimension of the minimal linear space it spans. Next,
they stratify low-genus curves according to the dimension of their
spans, as well as the types of their singularities, and
analyze the number of linear independent conditions imposed on curves
in each case. By controlling the dimension of each stratum inside the
appropriate mapping space of nondegenerate rational curves, they are
able to conclude that the projection $\Phi_d
\rightarrow M_d$ is never dominant over the locus of singular curves.

Johnsen and Kleiman also show that no reducible curves lie on the
general quintic. To do so, they bound the length of the
intersections of curve components, and thereby obtain lower bounds for
the codimension of such reducible
curves inside the Hilbert scheme. Those bounds allow them to conclude.

The next case to check, which is the subject of this work, is that of rational curves of degree $10$. It is
interesting in its relation to mirror symmetry. In fact (see
\cite[p.206]{Cox}), the finiteness of the Hilbert scheme of degree-$10$
rational curves on the general quintic smooth implies that the
instanton number $n_{10}$ is given by
\vspace{-.2cm}
\begin{equation}
n_{10}= 6 \times 17,601,000 +\#\{\text{smooth rational curves of degree
  $10$ in }F\}.\notag 
\end{equation}

A couple of comments are in
order. First of all, Clemens' original conjecture predicted that only
{\it smooth} rational curves lie on a general quintic $F$. Clemens'
conjecture in its revised form predicts that the six-nodal plane
quintics are the only singular, reduced, and irreducible rational
curves on $F$. On the other hand, mirror symmetry
includes Vainsencher's singular quintics in its count of rational curves of
degree five, but fails to count {\it six double covers} corresponding
to each of these. So we see both that $10$ is the first $d$ for which the
instanton number $n_d$ {\it fails to count smooth rational curves of
  degree $d$ on the general quintic $F$}, and that the discrepancy between $n_{10}$ and the actual number of smooth
rational curves is accounted for by double covers of nodal
plane quintics.

In this paper, we extend the results of Johnsen and Kleiman to include
curves of degree $10$. In doing so we adopt their basic strategy to show
that 
the only
rational, reduced, irreducible curves of degree $10$ on a general
quintic hypersurface $F$
are smooth, and that there are only finitely many such smooth curves. Since
the results of \cite{GLP} are less useful than before, however, we are forced to use a different
technique to estimate the dimensions
of the fibres of the projection $I_{10} \rightarrow M_{10}$. In order to make our dimension
estimates, we study the {\it generic initial ideals} which result from degenerations of rational curves.
By a combinatorial analysis, we are able to control the cohomology of generic initial ideals, and
deduce uniform bounds on $h^0(\mathcal{I}_{C/\mathbb{P}^4}(5))$, where $C$ is a rational
curve of degree $10$.

There are several reasons why studying the combinatorics of the generic initial ideal, or gin, is useful
(and feasible). First, it is a monomial ideal fixed by the
action of the group of upper triangular matrices, so it has a rather simple
structure. In addition, its Castelnuovo--Mumford regularity is the same
as that of (the ideal sheaf of) the curve $C$. Finally, the
presentation of the generic initial ideal $\mbox{gin }\mathcal{I}_{C}$ for the
reverse-lexicographic term order is closely related to the
presentation of $\mbox{gin }\mathcal{I}_{C \cap H/H}$, where $H$ is a general hyperplane. Since the general hyperplane
section of a curve $C$ is a collection of points in uniform position
(see, for example, \cite[p. 85]{Harris}),
there are natural restrictions on the Hilbert functions and regularity
properties of the corresponding ideals. 

By systematically using all of the preceding
considerations, we are able to obtain relatively robust bounds on the higher cohomology of $\mathcal{I}_{C/\mathbb{P}^4}(5)$, when $C$ is the image of any map $f\:\mathbb{P}^1 \rightarrow \mathbb{P}^4$ outside
a certain special locus inside the mapping space $M_{10}$. By combining
our bounds on cohomology with the estimates on parameter-space dimensions given in
\cite{Klei}, we are able to extend the results of \cite{Klei} to
degree $10$. To our knowledge, our method constitutes a new application of generic initial ideals
to problem-solving in algebraic geometry; hopefully, it will be a
useful tool for work on other problems, too. 

Hereafter, a ``generic'' quintic hypersurface in $\mathbb{P}^4$
means one cut out by a polynomial $F$ belonging to a Zariski-dense open
subset of the family $\mathbb{P}^{125}$ of degree-5 polynomials on
$\mathbb{P}^4$. We always work over the complex numbers
$\mathbb{C}$. Moreover, the abbreviations $\mathcal{I}_C$
and $\mathcal{I}_{\Gamma}$ are used in place of
$\mathcal{I}_{C/\mathbb{P}^4}$ and $I_{\Gamma/\mathbb{P}^3}$, respectively, for ideal sheaves of curves $C
\subset \mathbb{P}^4$ and their hyperplane sections. When we refer to the regularity of a variety
$X$, we mean the Castelnuovo--Mumford regularity of
the ideal sheaf,
$\mbox{reg }X:=\mbox{reg }\mathcal{I}_X$. Similarly, whenever $I$ is a
homogeneous ideal we take $H^i(I)$ and $\mbox{reg }I$ to mean the $i$th
cohomology group and the regularity, respectively, of the sheaf
associated to
$I$. The notation $I_{X,m}$ denotes the $m$th graded piece of the
homogeneous ideal $I_X$. Similarly, $\mathcal{I}_X|_H:= \mathcal{I}_X \otimes \mathcal{O}_H$ denotes the
restriction of $\mathcal{I}_X$ to the hyperplane $H$. Likewise,
$(I_X)_{x_n}$ denotes $\cup_{i=1}^{\infty} I_X:(x_n^i)$, the
saturation of $I_X$ with respect to the element $x_n$.

A {\it minimal generator} of a
monomial ideal $I$ means a monomial that is minimal for
the partial order defined by divisibility. (A basic fact from
commutative algebra is that every monomial ideal
$I$ admits a unique set of minimal generators, and these generate
$I$.) The term {\it genus} $g(C)$ of a one-dimensional scheme $C$
always denotes arithmetic genus, that is, $1-\chi(C)$, where $\chi$ is
the Euler characteristic of $\mathcal{O}_C$. The
{\it degree} or {\it genus} of a homogeneous ideal $I$ refers to the
degree or genus of the
projective scheme defined by $I$.
Finally, {\it we only consider initial ideals for the
reverse-lexicographic, or revlex, term order}, so henceforth $\mbox{gin }I$
denotes the generic initial ideal of $I$ with respect to
reverse-lexicographic order. Revlex is distinguished among term orders
for being ``best-behaved'' with regard to intersections with
hyperplanes. (For a discussion of term orders, see
\cite[Ch. 15]{Eisenbud} and \cite{Green}.)

The paper is in three sections. In section~\ref{sec-one}, we recall
the basic theory of generic initial ideals, in order to show how to compute the basic
cohomological invariants of such ideals. To bound $h^1(\mathcal{I}_C(5))$, we
first describe the initial ideals of general hyperplane sections of
irreducible curves of degree $10$, via a result of
\cite{Ballico} that bounds their Castelnuovo--Mumford regularity. Next we use a
result of \cite{Bayer2} to relate
the initial ideals of curves to those of their hyperplane sections. 

In section~\ref{sec-two}, the strong form of Clemens' conjecture for irreducible
rational curves of degree $10$ in $\mathbb{P}^{4}$ is
proved. To do so, we establish bounds on
$h^1(\mathcal{I}_C(5))$ case by case for rational curves $C$
of degree $10$. To conclude that the bounds we obtain are adequate, we appeal to
a result of Verdier (see \cite[Thm., p.139]{Ve} and \cite[Thm.1, p.181]{Ra}) that describes the stratification
of the space of all morphisms $f\:\mathbb{P}^1 \rightarrow \mathbb{P}^4$ of
a given degree according to the splitting of the
restricted tangent bundle 
$f^*T_{\mathbb{P}^4}$.

Section~\ref{sec-three} extends
Johnsen and Kleiman's analysis of reduced, reducible, and connected rational curves to cover
those of degree $10$. As in \cite{Klei}, we stratify curves
according to the length of the intersection of their
components. To bound the dimensions of these loci, we use a simple
argument involving the relative Hilbert scheme of morphisms, in place
of the arguments in local coordinates given in \cite{Klei}. Regularity
arguments from the latter paper carry over, with the exception of a regularity lemma for unions of nodal
quintics. We prove that such curves are $6$-regular (so, in particular,
they verify $h^{1}(\mathcal{I}_C(5))=0$) by showing their initial ideals are generated by polynomials of
degree $6$ or less.

Many of the methods discussed in this paper carry over to
degree $11$; however, one needs to work even harder to bound $h^{1}\mathcal{I}_{C}(5)$. In a preprint to be
distributed shortly, the
 author verifies Clemens' conjecture in degree $11$, by showing that those curves which do not
 satisfy suitable bounds on $h^{1}\mathcal{I}_{C}(5)$ necessarily lie on a large number of linearly independent
 hypercubics in $\mathbb{P}^{4}$. A liaison-theoretic argument shows that such curves necessarily lie on
 surfaces of degree at most $8$, and the author concludes by showing that the locus of rational curves lying on such 
 surfaces has small dimension.
\section*{Acknowledgements.} I am grateful to I. Coskun, J. Harris, and
A. Iarrobino, for providing me with valuable
advice during the preparation of this paper. I am especially grateful to S. Kleiman,
for an unending supply of input, encouragement, and patience.

\section{Generic initial ideals}\label{sec-one}
Given any homogeneous ideal $I$ in $n+1$ variables $x_i$,
together with a
choice of
partial order $<$ on the monomials $m(x_i)$, the {\it initial ideal} of $I$ with respect to $<$ 
is, by definition, generated by leading terms of elements in $I$. So
for any subscheme $C$ of $\mathbb{P}^n$, upper-semicontinuity implies that the regularity and values $h^j$
of the ideal sheaf $\mathcal{I}_C$ are majorized by those of the
(sheaf associated to the) initial ideal
$\mbox{in}(\mathcal{I}_C)$ with respect to any partial order on the
monomials of $\mathbb{P}^n$. On the other hand, Macaulay's theorem establishes that the Hilbert
functions of $\mathcal{I}_C$ and $\mbox{in}(\mathcal{I}_C)$ agree. If, moreover,
we replace $C$ by
a {\it general} $\mbox{PGL}(n+1)$-translate of $C$, then
\begin{equation*}
\mbox{reg}(\mathcal{I}_C) = \mbox{reg}(\mbox{in}(\mathcal{I}_C)),
\end{equation*}
and $\mbox{in}(\mathcal{I}_C)$ is called the {\it generic initial
  ideal}
of $C$, written $\mbox{gin}(\mathcal{I}_C$). 

A theorem of Galligo establishes
  that there is a unique generic initial ideal associated to any
  ideal. Further, Bayer and Stillman showed that the following three statements are equivalent.
\begin{enumerate}
\item The monomial ideal $\mbox{gin }\mathcal{I}_C$ is saturated, in
  the ideal-theoretic sense. In other words, for all polynomials $g \in \mbox{gin
  }\mathcal{I}_C$, $m \in (x_0,\dots,x_n)$,
\vspace{-.2cm}
\begin{equation*}
g \cdot m \in  \mbox{gin
  }\mathcal{I}_C \Rightarrow g \in \mbox{gin
  }\mathcal{I}_C.
\vspace{-.25cm}
\end{equation*}
\item The homogeneous ideal $I_C$
  is saturated.
\item No minimal generator of $\mbox{gin
  }\mathcal{I}_C \subset \mathbb{C}[x_0,\dots,x_n]$ is divisible by $x_n$.
\end{enumerate}
For proofs, see \cite[Thm.~2.30]{Green}.
More generally, let $I \subset \mathbb{C}[x_0,\dots,x_n]$ be any ideal that
is {\it Borel-fixed}, i.e., fixed under the action of upper triangular
matrices $\mathcal{T} \subset \mbox{PGL}(n+1)$; then $I$ is
saturated if and only if none of its minimal generators is divisible by $x_n$.
Hereafter, we work exclusively
with {\it saturated} generic initial ideals. In the subsections to follow, we
reduce the problem of bounding $h^1(\mathcal{I}_C(5))$, for reduced curves
$C$, to the problem of bounding $h^1(\mathcal{I}_X(5))$, for $X$ belonging to
a certain finite set of Borel-fixed ideals. We also give an explicit procedure for computing
$h^1(\mathcal{I}_X(5))$.

\subsection{Hyperplane sections of nondegenerate irreducible curves} \label{section1.1}
In this subsection, we examine monomial ideals that arise
as gins of hyperplane sections, or {\it hyperplane gins} for short, of nondegenerate integral rational
curves $C$ of degree $10$. More to the point, we examine monomial ideals
corresponding to hyperplane gins that are {\it general} in a sense we
will now make precise.

Begin by fixing $C$. Choose general homogeneous coordinates
 $x_0,\dots,x_4$ for $\mathbb{P}^4$, and set
$H:=\{x_4=0\}$. Then $H$ is a general hyperplane, whose intersection
with $C$ is a set of ten points in uniform position: the set's
monodromy is the full symmetric group on ten letters. Associated to
the saturated ideal defining $C \cap H$, there is a corresponding
Borel-fixed monomial
ideal in $\mathbb{C}[x_0,\dots,x_3]$, the saturated generic initial
 ideal $\mbox{gin
}\mathcal{I}_{C \cap H}$. This ideal is invariant under general
linear transformations of the coordinates $x_0,\dots,x_3$. We therefore call it {\it the }hyperplane gin
of $C$. 

As explained in the paragraphs preceding this subsection, saturatedness and
Borel-fixity together imply that the hyperplane gin of $C$ is minimally
generated by
monomials in the first three variables $x_0$, $x_1$, and $x_2$. Moreover, because it cuts
out a zero-dimensional scheme, $\mbox{gin }\mathcal{I}_{C \cap H}$ has some minimal generator of the form
$x_2^{\lambda}$ with $\lambda>0$. Indeed, if $I:=\mbox{gin
}\mathcal{I}_{C \cap H}$ contained no such minimal generator, then the vanishing
locus $V(I)$ would contain a $1$-dimensional projective subspace of
vectors
$\{(0,0,0,a_2,a_3):a_1,a_2 \in \mathbb{C}\}$, which is absurd according to the dimension theorem
\cite[Thm.9.3.11]{CLO} for monomial ideals.

On the other hand, Borel-fixity implies the hyperplane gin
is combinatorially simple: Recall that an ideal $I$ is Borel-fixed if and only
if for every monomial $P$,
\begin{equation}
P^{*}:=x_i/x_j \cdot P \text{ belongs to }I \text{ whenever }i < j;
\notag \end{equation}
see \cite[Thm.15.23]{Eisenbud}.
In constructing minimal generating sets for Borel-fixed monomial
ideals, we make systematic use of this fact without mentioning
it.

Any set of minimal generators for the general
hyperplane gin of $C$ may be represented by a {\it tree} $T$, or directed
graph
without cycles. Namely, fix an alphabet
$\mathcal{A}:=\{\emptyset,x_0,x_1,x_2\}$, and consider the set of all
trees with root vertices $\emptyset$ and all other vertices labeled by either
$x_0,x_1,\text{ or }x_2$. Call
terminal vertices {\it leaves}. Then the unique path from the root vertex
$\emptyset$ to any leaf labeled $x_{i_l}$ determines a sequence
of vertices $\emptyset, x_{i_1},x_{i_2},\dots,x_{i_l}$.  (We will only be concerned
with trees associated to nonzero ideals, so $\emptyset$ will never
arise as a leaf.) It's natural to
interpret the string $x_{i_1}x_{i_2} \dotsm x_{i_l}$ as a polynomial of degree
$l$. Likewise, if the unique path from the root vertex to a given
vertex $v$ involves $d$ edges, then we say that {\it $v$ has degree
  $d$}. A {\it rewriting rule} applied to a leaf $v$ is then a formal operation on $T$ that
at $v$
glues a certain number of new edges $e_i(v)$ terminating in vertices $v_i(v)$. The result is a new tree $T^{\prime}$
in which the vertices $v_i(v)$ are leaves of degree $d+1$. Say that a vertex $V_1$ {\it dominates} a vertex $V_2$ if it
is closer to the root vertex $\emptyset$. We also stipulate that, if
$V_i$ with label $x_i$ dominates $V_j$ with label $x_j$, then $i \leq
j$. 

For each homogeneous monomial ideal $I$, the unique minimal generating
set of $I$
determines a tree $T(I)$ whose leaves correspond to minimal generators
of $I$. The assignment $I \mapsto T(I)$, moreover,
is unique. Therefore, we treat rewriting rules interchangeably as
operations on either ideals or trees.
 
To clarify our terminology, here are a couple of simple examples of zero-dimensional subschemes of $\mathbb{P}^3$ and the
minimal generating sets of the corresponding ideals. There is exactly
one zero-dimensional Borel-fixed ideal of degree 1, namely
$I_1=(x_0,x_1,x_2)$. Similarly, $I_2=(x_0,x_1,x_2^2)$ is the unique zero-dimensional
Borel-fixed ideal of degree $2$ . Note that a
minimal generating set for $I_2$ may be obtained from one for $I_1$ by
exchanging the minimal generator $x_2$ of $I_1$ for the minimal
generator $x_2^2$ of $I_2$. In other words, a set of minimal
generators for $I_2$ may be obtained from a set of minimal generators
for $I_1$ by replacing $x_2$ with $x_2^2$ and keeping the other
generators. The corresponding operation on trees adds a new vertex
$v(x_2^2)$ to $T(I_1)$ and connects $v(x_2^2)$ with the vertex $v(x_2)$
corresponding to $x_2 \in I_1$ along a new edge. In other words, to
obtain $T(I_2)$ from $T(I_1)$ it suffices to apply a single rewriting
at the vertex $x_2$, that we denote by $x_2
\mapsto x_2^2$.

In general, as we show in the next subsection, any hyperplane section
gin $I$ of a curve may be constructed by applying a sequence of rewriting
rules of the form $X_k \mapsto Y_{k+1}$ to the ``empty ideal,''
$\emptyset$. Here $X_k$ is a single
monomial and $Y_{k+1}$ is a set of monomials that ``replace'' $X_k$. In other
words, $\emptyset$ and $I$ fit into a sequence of monomial
ideals
$I_0=\emptyset, \dots, I_n=I$ where $I_{k+1}$ is the ideal generated
by all monomial minimal generators of $I_{k}$ except $X_k$, together with a set of
monomials $Y_{k+1}$ of degree $\deg(X_k)+1$. 

\subsection{Rewriting rules for nondegenerate hyperplane gins}
\begin{table}[section]
\caption{$\Lambda$-rules for nondegenerate curves in $\mathbb{P}^4$}\label{Table $1$}
\vspace{-0.3cm}
$$\begin{array}{ll}
\hline
1.\rule{0 pt}{13 pt}&x_0^e
\mapsto (x_0^{e+1},x_0^ex_1,x_0^ex_2),\\
\vspace{.1cm}
2.&x_0^ex_1^f \mapsto (x_0^ex_1^{f+1},x_0^ex_1^fx_2), \\
\vspace{.1cm}
3.&x_0^ex_1^fx_2^g \mapsto
x_0^ex_1^fx_2^{g+1}, \text{and an {\it initial rule}}\\
4.&\emptyset \mapsto
(x_0,x_1,x_2).
\vspace{-.3cm} 
\end{array}$$
\end{table} 

Fix a nondegenerate curve $C$, together with a hyperplane $H$ that is
general in the sense that $C \cap H$ is a collection of 10 distinct
points in uniform position. Letting $I$ denote the corresponding hyperplane
gin, we see that $I$ is a homogeneous Borel-fixed ideal in the
coordinates $x_i, 0 \leq i \leq 3$, of $H$.
It follows from the discussion of subsection~{1.2} that $I$
is saturated if and only if none of its minimal generators is
divisible by $x_3$, and such a saturated ideal defines a
zero-dimensional scheme if and only if some minimal generator of $I$
is of the form $x_2^{\lambda}$. 

Our next task is to describe a set of rewriting rules for hyperplane
gins, or {\it $\Lambda$-rules}, that accounts for all
possible $I$ arising as above. A candidate for a complete list of $\Lambda$-rules is given in
Table~\ref{Table $1$}. Before showing that it is complete, a few words are in order.

First of all, it is clear that the proposed $\Lambda$-rules preserve
saturatedness and zero-dimensionality. Moreover, for every $\Lambda$-rule $X_k \mapsto Y_{k+1}$
that replaces the ideal $I_k$ by $I_{k+1}$,
the set of monomials $Y_{k+1}$ are in fact {\it minimal generators}
of $I_{k+1}$. Similarly, $\Lambda$-rules preserve
Borel-fixity because the combinatorial criterion that characterizes
the property is satisfied
at every step. So applying a sequence of $\Lambda$-rules to the set of
minimal generators of a hyperplane gin yields a set of minimal
generators for a Borel-fixed ideal of dimension zero.
 
As explained in the preceding subsection, it is also useful to interpret the $\Lambda$-rules
graphically. Namely, applying the $\Lambda$-rule $X_k \mapsto Y_{k+1}$ alters a tree of minimal generators
for $I_k$ by gluing $\mbox{Card}(Y_{k+1})$ new edges onto the leaf
corresponding to $X_k$. Therefore, at the risk of creating
additional confusion,
we will identify leaves $V$ with the monomial minimal
generators to which they correspond, i.e., by the unique sequence of vertices linking
$\emptyset$ to $V$. By doing so, we mean to make more transparent the correspondence
between gluing operations on trees of minimal generators and
replacement operations on sets of minimal generators.   

Finally we show that the $\Lambda$-rules comprise a complete set of rewriting rules for
hyperplane gins.
\newtheorem{Lemma}{Lemma}[subsection]
\begin{Lemma}
The minimal generating set of every nondegenerate hyperplane gin
may be realized by applying $\Lambda$-rules.
\end{Lemma}

\begin{proof}
Let $C$ and $\Gamma$ respectively denote a nondegenerate curve $C
\subset \mathbb{P}^4$ and a general
hyperplane section of the curve. Recall from subsection~\ref{1.2} that $\mbox{gin}(I_{\Gamma})$ is
Borel-fixed, and has a minimal generator of the form
$x_2^{\lambda}$, which is necessarily of maximal degree among all
minimal generators of $\mbox{gin}(I_{\Gamma})$.

Now consider the set of leaves $l$ of maximal degree $\lambda$ in $T(\mbox{gin}(I_{\Gamma}))$. Letting
$v_d(l)$ denote the vertex dominating $l$, note that $v_d(l)$ verifies
this combinatorial property: {\it if $v_d(l)$ is labeled $x_0$
(resp. $x_1,x_2$), then $v_d(l)$ dominates a $3$-tuple of vertices
  $x_0,x_1,x_2$ without omissions (resp., a pair $x_1,x_2$ or a
  singleton $x_2$)}. If $v_d(l)$ is labeled $x_2$, then the property
holds trivially, by the conventions for trees established in the previous subsection. In the other two cases, the observation follows from
  Borel-fixity, together with the fact that $x_2^{\lambda}$ lies in $\mbox{gin}(I_{\Gamma})$.

Note, moreover, that the tree obtained from
$T(\mbox{gin}(I_{\Gamma}))$ by ``pruning'', or removing all leaves of maximal degree and
contracting the corresponding edges to their dominating vertices, is
also Borel-fixed, saturated, and zero-dimensional. By induction on
degree, it follows that every vertex of $T(\mbox{gin}(I_{\Gamma}))$
verifies the combinatorial property of the previous paragraph.

On the other hand, any tree generated by $\Lambda$-rules may
be pruned, as pruning is clearly an inverse for rewriting.
The same argument by induction on degree invoked in the preceding paragraph now shows that those trees
obtained by $\Lambda$-rules are exactly those verifying the
combinatorial property above.
\end{proof} 

Now let $X$ be any subscheme of $\mathbb{P}^n$, and let
$\Gamma$ be a hyperplane section defined by $\Gamma:=X \cap H$
where $H$ is cut out by a general linear form. Make a general choice of coordinates $x_0,\dots,x_n$ on $\mathbb{P}^n$
with respect to which $H$ is defined by $x_n=0$. As explained in \cite[p. 163]{Green}, 
the surjectivity of the restriction map
\begin{equation*}
I_{X,m} \rightarrow I_{\Gamma,m},
\end{equation*}
for $m$ sufficiently large, implies that
\begin{equation*}
\mbox{sat}(I_{X}|_H)= I_{\Gamma}.
\end{equation*}
According to \cite[Prop. 2.21]{Green}, the latter equality of ideals may be formulated
as the following fact relating the generic initial ideal of $X$ to the generic initial ideal of
a general hyperplane $\Gamma$.
\newtheorem{Fact}{Fact}
\begin{Fact}\label{sat_lemma}
The
saturation of $\mbox{gin}(\mathcal{I}_X|_H)$ with respect to $x_{n-1}$
is equal to $\mbox{gin}(\mathcal{I}_{\Gamma})$.
\end{Fact}

Note that the transformation $\mbox{gin}(\mathcal{I}_X|_H) \mapsto
\mbox{gin}(\mathcal{I}_X|_H)_{x_{n-1}}$ induces a graphical
transformation $T(\mbox{gin}(\mathcal{I}_X|_H)) \mapsto
T(\mbox{gin}(\mathcal{I}_X|_H)_{x_{n-1}})$ given by iteratively
contracting all
edges in the tree dominated by vertices that dominate
$x_{n-1}$-labeled vertices, until no leaves labeled $x_{n-1}$
remain. On the other hand, the contractions taking
$T(\mbox{gin}(\mathcal{I}_X|_H))$ to
$T(\mbox{gin}(\mathcal{I}_X|_H)_{x_{n-1}})$ may be viewed as
transformations of trees of (minimal
generators of) homogeneous ideals of
$\mathbb{C}[x_0,\dots,x_n]$ that are inverse to certain
rewriting rules, which we call {\it $C$-rules} and display in Table~\ref{Table $2$}. Therefore, Fact~\ref{sat_lemma} may be reformulated in the
following useful way.

\begin{table}
\caption{$C$-rules for nondegenerate irreducible curves}\label{Table $2$}
\vspace{-.3cm}
$$\begin{array}{ll}
\hline
\vspace{0.1cm}
1.&x_0^e
\mapsto x_0^e \cdot(x_0,x_1,x_2,x_3),\\
\vspace{.1cm}
2.&x_0^ex_1^f \mapsto x_0^ex_1^f \cdot(x_1,x_2,x_3), \text{ and }\\
\vspace{.1cm}
3.& x_0^ex_1^fx_2^g \mapsto
x_0^ex_1^fx_2^g \cdot(x_2,x_3) \notag
\vspace{-.3cm}
\end{array}$$
\end{table}

{\it The tree corresponding to any
generic initial ideal defining a nondegenerate subscheme $X \subset \mathbb{P}^4$ of dimension 1 is
obtained from the $\mbox{gin}$ of a general hyperplane section of
$X$ by applying a sequence of $C$-rules.}

Since hyperplane gins are generated by $\Lambda$-rules, we have proved the following statement, which is
essential for what follows and will be used
hereafter without comment.

{\it Every Borel-fixed monomial ideal associated to a
curve $C \subset \mathbb{P}^4$
may be represented by a tree obtained from the empty symbol by applying a
certain sequence of $\Lambda$-rules, followed by a certain sequence of
$C$-rules.}

Likewise, we will make use of the following result of Bayer-Stillman (see
\cite[Thm.~2.27]{Green}):
{\it The regularity of any Borel-fixed homogeneous ideal $I$ is equal
  to
\begin{equation*}
\max_P\{\,\deg P\mid  P\text{ is a minimal generator of }\mbox{gin(I)}\,\}.
\end{equation*}
}
Here, as everywhere else in this paper, the underlying term order on monomials
is assumed to be revlex; indeed the proposition fails to hold for
arbitrary term orders.

The following lemma permits us to identify all Borel-fixed ideals
associated to general hyperplane sections of subschemes of
$\mathbb{P}^4$.  
\begin{Lemma}\label{zerodeg}
The number of nonleaf vertices in the tree defining a zero-dimensional
subscheme of $\mathbb{P}^3$ is equal to the scheme's degree. 
\end{Lemma}

\begin{proof}
We proceed by induction on the degree of the subscheme $\Gamma \subset
\mathbb{P}^3$. Without loss of
generality, we may assume
$\Gamma$ is nondegenerate, since the restriction $\Gamma|_L$ of $\Gamma$ to its
linear span $L$ satisfies $\deg(\Gamma|_L)=\deg(\Gamma)$. Thus, $\mbox{gin}(I_{\Gamma})$ is a monomial
ideal whose minimal generators are polynomials in $x_0,x_1,\text{ and }x_2$. Moreover, since $\Gamma$ is zero-dimensional,
$\mbox{gin}(I_{\Gamma})$ has one minimal generator of the form
$x_2^{\lambda}$, for some positive integer
$\lambda$. Furthermore, $\lambda$ is the maximal total degree of any
minimal generator of $\mbox{gin}(I_{\Gamma})$.

If $\deg(\Gamma)=1$, the claim is trivial, since the only
nondegenerate zero-dimensional subscheme of $\mathbb{P}^3$ of degree
one is the point, whose ideal (up to projective linear transformation) is $(x_0,x_1,x_2)$. If $\deg(\Gamma)>1$, then $\mbox{gin}(I_{\Gamma})$ may be obtained
from the gin of the point by applying finitely many $\Gamma$-rules, so
by induction it suffices to check the following statement:

{\it Let $\Gamma$ and $\Gamma^{\prime}$ be zero-dimensional subschemes
  of $\mathbb{P}^3$ defined by Borel-fixed ideals $I$ and $I^{\prime}$ of
  $\mathbb{C}[x_0,x_1,x_2,x_3]$, respectively. If $I^{\prime}$ is
  obtained from $I$ by applying a single
  $\Lambda$-rule to $I$, then $\deg(\Gamma^{\prime})=\deg(\Gamma)+1$.}

To check the latter claim, first note that by the regularity
criterion \cite[Thm.~2.27]{Green} of Bayer-Stillman cited above, the ideal sheaves $\mathcal{I}$ and
$\mathcal{I}^{\prime}$ associated to the two homogeneous ideals in
question are $\lambda^{\prime}$-regular, where 
\begin{equation*}
\lambda^{\prime}:=\max_P\{\,\deg(P)\mid P \text{ is a minimal generator of } I^{\prime}\,\}.
\end{equation*}
Therefore, there are exact sequences
on global sections:
\begin{equation*}
0 \rightarrow H^0(\mathcal{I}(\lambda^{\prime})) \rightarrow
H^0(\mathcal{O}_{\mathbb{P}^3}(\lambda^{\prime})) \rightarrow H^0(\mathcal{O}_{\Gamma}(\lambda^{\prime}))
\rightarrow 0
\end{equation*} and
\begin{equation*}
0 \rightarrow H^0(\mathcal{I}^{\prime}(\lambda^{\prime})) \rightarrow
H^0(\mathcal{O}_{\mathbb{P}^3}(\lambda^{\prime})) \rightarrow H^0(\mathcal{O}_{\Gamma^{\prime}}(\lambda^{\prime}))
\rightarrow 0.
\end{equation*}

Moreover, since $\Gamma$ and $\Gamma^{\prime}$ are zero-dimensional,
their degrees are equal to their Euler characteristics. Since
$\mathcal{I}(\lambda^{\prime}), \mathcal{I}^{\prime}(\lambda^{\prime})$, and
$\mathcal{O}_{\mathbb{P}^3}(\lambda^{\prime})$ have no higher cohomology,
$\chi=h^0$ for each of these sheaves.

Because $\chi$ is additive across short exact
sequences, it follows that
\begin{equation*}
\chi(\Gamma)=h^0(\mathcal{O}_{\mathbb{P}^3}(\lambda^{\prime}))-H^0(\mathcal{I}(\lambda^{\prime}))
\text{ and }
\chi(\Gamma^{\prime})=h^0(\mathcal{O}_{\mathbb{P}^3}(\lambda^{\prime}))-h^0(\mathcal{I}^{\prime}(\lambda^{\prime})).
\end{equation*}
But clearly
$h^0(\mathcal{I}^{\prime}(\lambda^{\prime}))=h^0(\mathcal{I}(\lambda^{\prime}))-1$,
since $I^{\prime}$ is a saturated homogeneous ideal with one fewer
monomial generator (not necessarily minimal!) in degree $\lambda^{\prime}$ than $I$, which is also
saturated and homogeneous. It follows immediately that
$\chi(\Gamma^{\prime})=\chi(\Gamma)+1$, whence $\deg(\Gamma^{\prime})=\deg(\Gamma)$+1.
\end{proof}

\begin{figure}
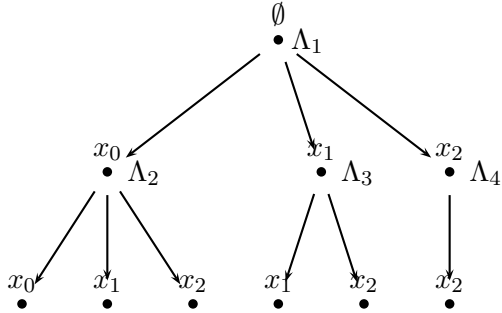

\begin{center}
\psset{tnpos=a,radius=2pt,nodesep=7pt,labelsep=2pt,levelsep=50pt,arrows=->}

\pstree{\Tdot*~{$\emptyset$}~[tnpos=r]{$\Lambda_1$}}{
\pstree{\Tdot*~{$x_0$}~[tnpos=r]{$\Lambda_2$}}{\Tdot*~{$x_0$}
\Tdot*~{$x_1$} \Tdot*~{$x_2$}}
\pstree{\Tdot*~{$x_1$}~[tnpos=r]{$\Lambda_3$}}{\Tdot*~{$x_1$} \Tdot*~{$x_2$}}
\pstree{\Tdot*~{$x_2$}~[tnpos=r]{$\Lambda_4$}}{\Tdot*~{$x_2$}}
}
\caption{The generic initial ideal of a rational normal
  quartic $C \subset \mathbb{P}^4$.} 
\label{figure: Figure 1}
\end{center}
\end{figure}

{\bf Example.} Figure~\ref{figure: Figure 1} shows the tree-representation of the hyperplane gin of a 
rational normal quartic in $\mathbb{P}^4=\mbox{Proj
}\mathbb{C}[x_0,x_1,x_2,x_3,x_4]$. Note that in this example, 
\begin{equation*}
\mbox{gin}(\mathcal{I}_C)=\mbox{gin}(\mathcal{I}_{C \cap H})^{\mbox{ext}},
\end{equation*}
where $H$ is any hyperplane generic with respect to $C$, and
$\mbox{gin}(\mathcal{I}_{C \cap H})^{\mbox{ext}}$ denotes the
extension of $\mbox{gin}(\mathcal{I}_{C \cap H})$ to $\mathbb{C}[x_0,\dots,x_4]$. The $\Lambda$-rules are
  marked as $\Lambda_i$, for $i=1,\dots,4$. Exactly four vertices are not leaves.

In this example the hyperplane gin is
\begin{equation*}
I=(x_0^2,x_0x_1,x_0x_2,x_1^2,x_1x_2,x_2^2),
\end{equation*}
the gin of four points in general position in
$\mathbb{P}^3$. These points arise as the generic hyperplane section
of a rational normal quartic, so according to Fact \ref{sat_lemma} the gin
of some rational normal quartic must be obtainable by the application
of some nonnegative number $m$ of $\Lambda$-rules from
$I$. Indeed, one can check (using, e.g., the computer algebra system
Macaulay2 of \cite{Ma}) that the gin of the rational normal quartic is
$(x_0^2,x_0x_1,x_0x_2,x_1^2,x_1x_2,x_2^2)$; so in this case, $m=0$.

To construct the tree in the example, we may proceed as
follows. First, apply the initial $\Lambda$-rule to obtain the ideal $(x_0,x_1,x_2)$. Next,
apply a rewriting rule at the vertex corresponding to the
generator $x_2$, replacing $x_2$ with $x_2^2$. The corresponding rule is $x_2
\mapsto x_2^2$. Similarly, replace $x_1$ with
$(x_1^2,x_1x_2)$ via $x_1 \mapsto (x_1^2,x_1x_2)$ and replace $x_0$ with
$(x_0^2,x_0x_1,x_0x_2)$ via $x_0 \mapsto (x_0^2,x_0x_1,x_0x_2)$.

\subsection{Hyperplane gins for nondegenerate curves}
In this subsection, $I$ denotes a saturated Borel-fixed ideal defining the gin of some general hyperplane section of some (fixed)
nondegenerate, irreducible degree-$10$
rational curve $C$. Throughout, we use $\mbox{Borel}(J)$ to denote the smallest
Borel-fixed ideal containing the ideal $J$. We shall see that the minimal generating set of
$I$ is subject to significant numerical restrictions, which together
imply there are very few possibilities
for $I$. Four out of the five corresponding possible trees, which we will implicitly
refer to in the course of establishing bounds on $h^1(\mathcal{I}_C(5))$, are given in figures 2 through 5. 

To begin, note that every general hyperplane
section $\Gamma$ of $C$ is (the reduced scheme associated
to) a set of ten points in
uniform position in a three-dimensional projective space. By a result of
\cite{Ballico} bounding the regularity of ideals of points in uniform position, $\mathcal{I}_{\Gamma}$ is
$4$-regular, so $I$ is minimally generated in degrees at most $4$.

{\it If $I$ has no quadratic generators}, then 
\begin{equation}
\begin{split}
I&=\mbox{Borel}(x_2^3)\\
&=(x_0^3,x_0^2x_1,x_0^2x_2,x_0x_1^2,x_0x_1x_2,x_0x_2^2,x_1^3,x_1^2x_2,x_1x_2^2,x_2^3),
 \notag
\end{split}
\end{equation}
by saturatedness,
Borel-fixity, and Lemma \ref{zerodeg}. See Figure~\ref{figure: Figure $2a$}.

\begin{figure}[htbp]
\begin{center}
\psset{tnpos=a,nodesep=7pt,labelsep=1pt,levelsep=40pt,arrows=->}

\pstree{\Tdot*~{$\emptyset$}}{

\pstree{\Tdot*~{$x_0$}}{ 

\pstree{\Tdot*~{$x_0$}} 
{\Tdot*~{$x_0$} \Tdot*~{$x_1$} \Tdot*~{$x_2$}}

\pstree{\Tdot*~{$x_1$}} 
{\Tdot*~{$x_1$} \Tdot*~{$x_2$}}

\pstree{\Tdot*~{$x_2$}} 
{\Tdot*~{$x_2$}}
}

\pstree{\Tdot*~{$x_1$}}{ 

\pstree{\Tdot*~{$x_1$}} 
{\Tdot*~{$x_1$} \Tdot*~{$x_2$}}

\pstree{\Tdot*~{$x_2$}} 
{\Tdot*~{$x_2$}}
}

\pstree{\Tdot*~{$x_2$}}{ 

\pstree{\Tdot*~{$x_2$}} 
{\Tdot*~{$x_2$}}
}
}

\caption{$\mbox{Borel}(x_2^3)$}
\label{figure: Figure $2a$}

\end{center}
\end{figure}

Similarly, {\it if $I$ has exactly one quadratic generator}, then
\begin{equation}
I=\mbox{Borel}(x_1x_2^2)+(x_0^2,x_2^4).
\notag \end{equation}
See Figure~\ref{figure: Figure $3$}.

\begin{figure}[htbp]
\begin{center}
\psset{tnpos=a,nodesep=7pt,labelsep=1pt,levelsep=30pt,arrows=->}

\pstree{\Tdot*~{$\emptyset$}}{

\pstree{\Tdot*~{$x_0$}}{ 
\Tdot*~{$x_0$}
\pstree{\Tdot*~{$x_1$}}{\Tdot*~{$x_1$} \Tdot*~{$x_2$}}
\pstree{\Tdot*~{$x_2$}}{\Tdot*~{$x_2$}}}

\pstree{\Tdot*~{$x_1$}}{ 

\pstree{\Tdot*~{$x_1$}} 
{\Tdot*~{$x_1$} \Tdot*~{$x_2$}}

\pstree{\Tdot*~{$x_2$}} 
{\Tdot*~{$x_2$}}
}

\pstree{\Tdot*~{$x_2$}}{ 

\pstree{\Tdot*~{$x_2$}}{ 

\pstree{\Tdot*~{$x_2$}}{ 
\Tdot*~{$x_2$}}
}
}
}
\caption{$\mbox{Borel}(x_1x_2^2)+(x_0^2,x_2^4)$}
\label{figure: Figure $3$}

\end{center}
\end{figure}

\begin{figure}[htbp]
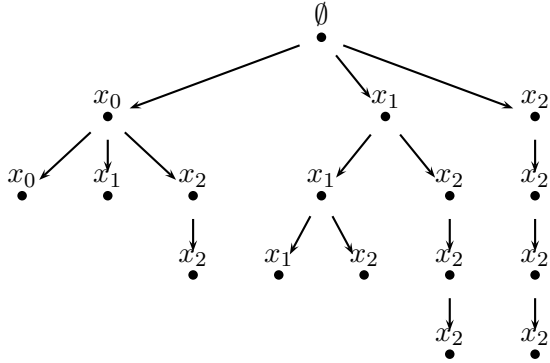

\begin{center}
\psset{tnpos=a,nodesep=7pt,labelsep=1pt,levelsep=30pt,arrows=->}

\pstree{\Tdot*~{$\emptyset$}}{

\pstree{\Tdot*~{$x_0$}}{ 
\Tdot*~{$x_0$} \Tdot*~{$x_1$}
\pstree{\Tdot*~{$x_2$}}{\Tdot*~{$x_2$}}}

\pstree{\Tdot*~{$x_1$}}{ 

\pstree{\Tdot*~{$x_1$}} 
{\Tdot*~{$x_1$} \Tdot*~{$x_2$}}

\pstree{\Tdot*~{$x_2$}}{ 
\pstree{\Tdot*~{$x_2$}}{ 
\Tdot*~{$x_2$}}}}

\pstree{\Tdot*~{$x_2$}}{ 

\pstree{\Tdot*~{$x_2$}}{ 

\pstree{\Tdot*~{$x_2$}}{ 
\Tdot*~{$x_2$}}
}
}
}
\caption{$\mbox{Borel}(x_0x_1)+(x_1^3,x_1^2x_2,x_0x_2^2)+(x_1x_2^3,x_2^4)$}
\label{figure: Figure $4$}

\end{center}
\end{figure}

Next say that $I$ has exactly two quadratic generators. Borel-fixity
implies these are necessarily $x_0^2$ and $x_0x_1$.
Note that $I \neq \mbox{Borel}(x_2^3)+(x_0x_1,x_0^2)$ since the
corresponding tree has eight nonterminal vertices. Similarly, $x_1x_2^2$ cannot belong to $I$,
since otherwise degree considerations force
\begin{equation*}
I=\mbox{Borel}(x_0x_1)+\mbox{Borel
  }(x_1x_2^2)+(x_2^5). 
\end{equation*}
The latter ideal is not 4-regular, and therefore violates the main
result of \cite{Ballico}. Therefore,
{\it if $I$ has exactly two quadratic generators, then}
\begin{equation*}
I= \mbox{Borel}(x_0x_1)+(x_1^3,x_1^2x_2,x_0x_2^2)+(x_1x_2^3,x_2^4).
\end{equation*}
See Figure~\ref{figure: Figure $4$}.

On the other hand,{\it if $I$ has exactly three quadratic generators, }then these are
  either
$\{x_0^2,x_0x_1,x_0x_2\} \text{ or }
\{x_0^2,x_0x_1,x_1^2\}$, so
degree considerations, saturatedness, and Borel-fixity force
\begin{equation*}
I=\mbox{Borel}(x_0x_2)+\mbox{Borel}(x_2^4)+(x_1^3) \text{ or
  }I=\mbox{Borel}(x_0x_1)+(x_0x_2^3,x_1^2,x_1x_2^3,x_2^4).
\end{equation*}

\begin{figure}[htbp]
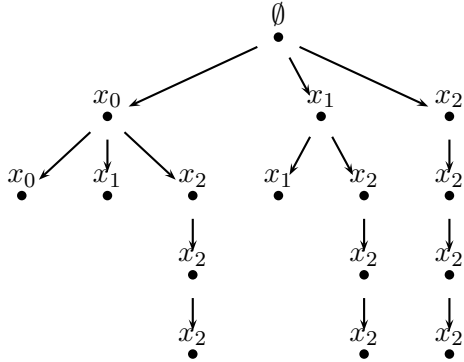

\begin{center}
\psset{tnpos=a,nodesep=7pt,labelsep=1pt,levelsep=30pt,arrows=->}

\pstree{\Tdot*~{$\emptyset$}}{

\pstree{\Tdot*~{$x_0$}}{ 
\Tdot*~{$x_0$} \Tdot*~{$x_1$} 
\pstree{\Tdot*~{$x_2$}}{
\pstree{\Tdot*~{$x_2$}}{\Tdot*~{$x_2$}}
}}

\pstree{\Tdot*~{$x_1$}}{ 
\Tdot*~{$x_1$} 
\pstree{\Tdot*~{$x_2$}}{ 
\pstree{\Tdot*~{$x_2$}}{
\Tdot*~{$x_2$}}}}

\pstree{\Tdot*~{$x_2$}}{ 

\pstree{\Tdot*~{$x_2$}}{ 

\pstree{\Tdot*~{$x_2$}}{ 
\Tdot*~{$x_2$}}
}
}
}
\caption{$\mbox{Borel}(x_0x_1)+(x_0x_2^3,x_1^2,x_1x_2^3,x_2^4)$}
\label{figure: Figure $5$}

\end{center}
\end{figure}

See Figure~\ref{figure: Figure $5$} for a presentation
  of $I=\mbox{Borel}(x_0x_1)+(x_0x_2^3,x_1^2,x_1x_2^3,x_2^4)$.

Note that {\it no zero-dimensional, saturated, and nondegenerate Borel-fixed ideal $I
  \subset \mathbb{C}[x_0,\dots,x_3]$ of degree
$10$ having at least four quadratic generators is 4-regular}. For, by
  Borel-fixity, if $I$ has at least four quadratic generators then
  these necessarily include
\begin{equation*}
x_0^2,x_0x_1,x_0x_2,\text{ and }x_1^2,
\end{equation*}
so that
\begin{equation*}
I=(x_0^2,x_0x_1,x_0x_2,x_1^2)+(x_1x_2^e)+(x_2^f),
\end{equation*}
for some $e$ and $f$. Borel-fixity implies $f \geq e$, but then
Lemma \ref{zerodeg} implies $e+f=10$, so $I$ has minimal generators in
degree $5$ or higher.

In sum, we have proved the following result.
\newtheorem{proposition}{Proposition}[subsection]
\begin{proposition}
$I$ has at most three quadratic generators, and is one
of the five monomial ideals listed above.
\end{proposition}

\subsection{Computing $h^1(\mathcal{I}_C(5))$ from a tree of minimal generators}
In this subsection we prove a few technical lemmas that will be used, in
the proof of Theorem~\ref{1.1}, to obtain bounds on
$h^1(\mathcal{I}_C(5))$.

\begin{Lemma}\label{L2}Let $C \subset \mathbb{P}^4$ denote any nondegenerate
  degree-$10$ integral curve. Then
  $h^1(\mathcal{O}_C(5))=0$.
\end{Lemma}

\begin{proof}
The Castelnuovo--Halphen genus bound (see, e.g.,
\cite[(1.1),p.27]{Ciliberto}) implies $g(C) \leq 9$, whence
$\deg(K_C)=2g(C)-2 \leq 16$. Since $C$ is integral (and in particular Cohen-Macaulay), Grothendieck
duality (see \cite{AK}) implies
\begin{equation}
h^1(\mathcal{O}_C(5))=h^0(K_C(-5)), \notag
\end{equation}
where $K_C$ denotes the
dualizing sheaf $C$. It therefore suffices
to show that $K_C(-5)$ has
no global sections.

Suppose $K_C(-5)$ has a section. Then there is an injection
\begin{equation*}
0 \rightarrow \mathcal{O}_C \rightarrow K_C(-5),
\end{equation*}
or equivalently, an injection of invertible sheaves
\begin{equation}\label{Kseq}
0 \rightarrow \mathcal{O}_C(5) \rightarrow K_C.
\end{equation}
In particular, we have 
\begin{equation*}
\chi(K_C) =\chi(\mathcal{O}_C(5))+ \chi(\mathcal{Q}),
\end{equation*}
where $\mathcal{Q}$ denotes the quotient of \eqref{Kseq}. Since \eqref{Kseq}
is generically an isomorphism, $\mathcal{Q}$ is supported at finitely
many points, so $\chi(\mathcal{Q})=h^0(\mathcal{Q})$ is nonnegative,
and therefore
\begin{equation}\label{chi_inequality}
\chi(K_C) \geq \chi(\mathcal{O}_C(5)).
\end{equation}
The right side of \eqref{chi_inequality} equals $5 \times 10 -g(C)+1$ and is therefore at
least $42$, whereas the left equals $2g(C)-2-g(C)+1$ and is therefore
at most $8$. So $h^0(K_C(-5))=0$ after all, and the lemma is proved.
\end{proof}

As explained in Section~\ref{section1.1}, it follows that
\begin{equation}
h^0(\mathcal{I}_C(5))=126-(5 \times 10-g+1)+i. \label{E1.3}
\end{equation}
On the other hand, $\mathcal{I}_C$ and $\mbox{gin }\mathcal{I}_C$ have the same
regularity, so since 
\begin{equation*}
h^2(\mathcal{I}_C(5))=h^1(\mathcal{O}_C(5))=0,
\end{equation*}
we also have
$h^2(\mbox{gin }\mathcal{I}_C(5))=0$. So in fact \eqref{E1.3} remains
true when $\mathcal{I}_C$ is replaced by $\mbox{gin}(\mathcal{I}_C)$. 

\begin{Lemma}\label{L3}Let $I \subset \mathbb{C}[x_0,\dots,x_4]$ be a Borel-fixed monomial ideal
  of dimension $1$, and $I^{\prime}$ the ideal obtained from $I$ by applying a
  single $C$-rule. Then $g(I^{\prime})=g(I)-1$.\end{Lemma}
(Recall from the introduction that $g(I)$ and $g(I^{\prime})$ denote the genera of
  the subschemes of $\mathbb{P}^4$ defined by $I$ and $I^{\prime}$.)

\begin{proof} Let $\mathcal{I}$
and $\mathcal{I}^{\prime}$ denote the sheaves associated to $I$ and
$I^{\prime}$, and let $C$ and $C^{\prime}$ denote the corresponding subschemes
of $\mathbb{P}^4$. Let $m$ be any positive integer for which all four
sheaves are $(m-1)$-regular. Then
$h^1(\mathcal{I}(m))=h^1(\mathcal{I}^{\prime}(m))=0$, so there are exact sequences of sections
\begin{equation*}
0 \rightarrow H^0(\mathcal{I}(m)) \rightarrow
H^0(\mathcal{O}_{\mathbb{P}^4}(m)) \rightarrow H^0(\mathcal{O}_C(m))
\rightarrow 0
\end{equation*} and
\begin{equation*}
0 \rightarrow H^0(\mathcal{I}^{\prime}(m)) \rightarrow
H^0(\mathcal{O}_{\mathbb{P}^4}(m)) \rightarrow H^0(\mathcal{O}_{C^{\prime}}(m))
\rightarrow 0,
\end{equation*}
Moreover, since
$h^1(\mathcal{O}_C(m))=h^1(\mathcal{O}_{C^{\prime}}(m))=0$, we have 
\begin{equation}\label{E1}
h^0(\mathcal{I}(m))=\binom{m+4}{m}-(10m-g(C)+1)
\end{equation} and
\begin{equation}\label{E2}
h^0(\mathcal{I}^{\prime}(m))=\binom{m+4}{m}-(10m-g(C^{\prime})+1)
\end{equation}
by the Riemann-Roch formula.
On the other hand, $\mathcal{I}$ and $\mathcal{I}^{\prime}$ are both saturated, so
$h^0(\mathcal{I}(m))$ and $h^0(\mathcal{I}^{\prime}(m))$ are exactly the numbers of
linearly independent monomials in the homogeneous ideals $I$ and $I^{\prime}$, respectively, of total degree
at most $m$. Hence
\begin{equation*}
h^0(\mathcal{I}^{\prime}(m))=h^0(\mathcal{I}(m))-1.
\end{equation*}
The
lemma now follows immediately from equations~\eqref{E1} and \eqref{E2}.
\end{proof}

Now let $C \subset \mathbb{P}^4$ be a nondegenerate degree-$10$
integral curve, and set $i:=h^1(\mathcal{I}_C(5))$.
\begin{Lemma}\label{4} The number of vertices in the tree representing $\mbox{gin
}\mathcal{I}_C(5)$ dominating vertices of degree greater than
  $6$ (equivalently, the number of rewritings applied to vertices of
  degree $6$ or greater) equals $i$. 
\end{Lemma}

\begin{proof} Let $v$ denote the number of rewritings applied to
  vertices of degree $6$ or greater. Recall (\cite[Thm.~2.27]{Green}) that any generic initial ideal $I$
for the reverse-lexicographic order is $m$-regular if and only if it
is minimally generated in degrees $m$ or less. So if $v=0$, then $C$
is $6$-regular, and therefore $i=0=v$.  Now say $v>0$. By
Lemma~\ref{L3}, each successive $C$-rewriting $r_i:T(I_i) \rightarrow
T(I_{i+1})$ satisfies
\begin{equation*}
g(C_{i+1})=g(C_i)-1, 
\end{equation*}
where $C_i$ and $C_{i+1}$ denote the subschemes of $\mathbb{P}^4$
defined by $I_i$ and $I_{i+1}$, respectively.
Now let $\deg(r_i)$ denote the degree of the vertex of $T(I_i)$ at which the rewriting $r_i$ is
applied. Note that if
$\deg(r_i) \geq 6$, then the number of linearly independent quintic
polynomials in $I_{i+1}$ equals the corresponding number in $I_i$. So, because rewritings preserve saturatedness,
\begin{equation*}
h^0(\mathcal{I}_{C_{i+1}}(5))=h^0(\mathcal{I}_{C_i}(5)).
\end{equation*} 
By \eqref{E1.3}, it
follows that if $r_i \geq 6$, then
\begin{equation*}
h^1(\mathcal{I}_{C_{i+1}}(5))=h^1(\mathcal{I}_{C_i}(5))+1.
\end{equation*} 
The lemma follows immediately by induction on $v$.
\end{proof} 
Note that by Lemma~\ref{L3}, $i$ also measures the failure of $C$ to impose linear
  independent conditions on quintic hypersurfaces.

\subsection{Nondegenerate curve gins associated to hyperplane gins}
In this subsection, we obtain restrictions on
minimal generating sets of
generic initial ideals of irreducible, nondegenerate rational curves of
degree $10$ in $\mathbb{P}^4$. To see how this is possible, fix a
choice $C$ of such a curve, with
$\Gamma$ a hyperplane section of $C$ defined by 
\begin{equation*}
\Gamma:=C \cap H, 
\end{equation*}
where $H$
is a general linear form. Make a general choice of coordinates $x_0,
\dots, x_4$ on $\mathbb{P}^4$ with respect to which $H$ is defined by
$x_4=0$. Let $I:= \mbox{gin}(\mathcal{I}_{\Gamma})$. Then by Fact \ref{sat_lemma}, $\mbox{gin
}\mathcal{I}_C$ is obtained via a sequence of
$C$-rewriting rules from the extension of $I \subset \mathbb{C}[x_0,\dots,x_3]$, to
$\mathbb{C}[x_0,\dots,x_4]$. 

As usual, let $i:=h^1(\mathcal{I}_C(5))$. The result of \cite[Thm.~3.1, p. 501]{GLP} implies that the ideal sheaf $\mathcal{I}_{C}$ is
$7$-regular unless $C$ admits an $8$-secant line (here points along the
line are counted with multiplicity). So {\it most} of the time,
\cite[Thm.~2.27]{Green} implies that $\mbox{gin
}\mathcal{I}_C$ is minimally generated by polynomials of degree at
most $7$. This fact limits the number of $C$-rewritings
that occur in degrees at least $6$ used to obtain $\mbox{gin
}\mathcal{I}_C$ from $I$, which in turn measures $i$, according to Lemma \ref{4}.

Let $C_{\Gamma}$ denote the cone with vertex
$(0,0,0,0,1)$ over the zero-dimensional scheme
defined by the vanishing of $I$ in $H$. Thus $C_{\Gamma}$ is one-dimensional, and its
minimal generators are exactly those of $I$; that is, $I_{C_{\Gamma}}$ is the
extension of $I$ to $\mathbb{C}[x_0,\dots,x_4]$. Let
$g_{\Gamma}:=g(C_{\Gamma})$. In the proof of Theorem~\ref{1.1}, we
will repeatedly use the following two technical results.
\begin{Lemma}\label{g+i}
$g_{\Gamma}-g$ is the number of $C$-rewritings applied to yield $\mbox{gin}(\mathcal{I}_C)$ from $I$, and $g+i \leq g_{\Gamma}$.
\end{Lemma}
\begin{proof}
The first statement is an immediate consequence of Lemma \ref{L3}. The
second statement therefore follows immediately whenever $i=0$. Note
that $\mbox{gin
}\mathcal{I}_C$ may be obtained in two steps:
\begin{enumerate}
\item Perform a number $r_1$ of
rewritings in degrees less than six.
\item Perform a number $r_2$ of
rewritings in degrees six or greater.
\end{enumerate}
Let $\stackrel{\sim}{C}$ denote
the scheme defined by the ideal that is the outcome of step
$1$. Clearly, we have 
\begin{equation*}
g_{\Gamma}-g(\stackrel{\sim}{C})=r_1
\end{equation*}
and
\begin{equation*}
g(\stackrel{\sim}{C})-g=r_2.
\end{equation*}
But Lemma \ref{4} implies that
\begin{equation*}
i=r_2,
\end{equation*}
so the second statement of the Lemma follows immediately. 
\end{proof}

Now, let $m$ be any positive integer such that $C$ is $m$-regular.
\begin{Lemma}\label{g_Gamma}
$g_{\Gamma}=10m+1-\binom{m+4}{4}+h^0(\mathcal{I}_{C_{\Gamma}}(m))$.
\end{Lemma}
\begin{proof}
Because $C_{\Gamma}$ is $m$-regular, the long exact sequence in cohomology associated to
\begin{equation*}
0 \rightarrow \mathcal{I}_{C_{\Gamma}}(m) \rightarrow
\mathcal{O}_{\mathbb{P}^4}(m) \rightarrow
\mathcal{O}_{C_{\Gamma}}(m) \rightarrow 0
\end{equation*}
shows that $\mathcal{O}_{C_{\Gamma}}$ is $m$-regular, too. Therefore,
there is an exact sequence
\begin{equation*}
0 \rightarrow H^0(\mathcal{I}_{C_{\Gamma}}(m)) \rightarrow
H^0(\mathcal{O}_{\mathbb{P}^4}(m)) \rightarrow
H^0(\mathcal{O}_{C_{\Gamma}}(m)) \rightarrow 0,
\end{equation*}
and $h^1(\mathcal{O}_{C_{\Gamma}}(m))=0$. It follows that
\begin{equation}\label{eqn1}
h^0(\mathcal{O}_{C_{\Gamma}}(m))=10 \times m-g_{\Gamma}+1,
\end{equation}
and also
\begin{equation}\label{eqn2}
h^0(\mathcal{O}_{C_{\Gamma}}(m))=\binom{m+4}{4}-h^0(\mathcal{I}_{C_{\Gamma}}(m)).
\end{equation}
To conclude, simply compare \eqref{eqn1} and \eqref{eqn2}.
\end{proof}

\section{Irreducible rational curves of degree $10$ in $\mathbb{P}^4$}\label{sec-two}
In this section, we'll prove the following theorem, which extends to
degree $10$ an
earlier result (\cite[Thm.3.1]{Klei}) of Johnsen and Kleiman's.

\newtheorem{thm}{Theorem}[section]
\begin{thm}\label{1.1}
The incidence
scheme $\Phi_d$ of smooth rational curves of degree $d$ at most $10$ on quintic
hypersurfaces $F \subset \mathbb{P}^4$ is irreducible. Moreover, any smooth curve $C$ lying on
a general quintic $F$ is embedded with normal bundle $\mathcal{O}_C(-1) \oplus
\mathcal{O}_C(-1)$. Furthermore, there are no rational and
singular, reduced, and irreducible curves of degree at most $10$ lying on
a general quintic in $\mathbb{P}^4$, other than
the six-nodal plane quintics.
\end{thm}

As noted in the introduction, there is an interesting corollary with
significance for mirror symmetry.

\newtheorem{cor}{Corollary}[section]
\begin{cor}\label{1.1cor}
The instanton number $n_{10}$ for a general quintic
threefold $F$ does not equal the number of smooth rational curves of
degree $10$ that lie on $F$.
\end{cor}

To prove the theorem, we'll borrow extensively from Johnsen and
Klei\-man's work; they prove the theorem in degree at most $9$.

\subsection{Initial reductions}
Let $M_{10}^{r,g}$ denote the affine space of parameterized
mappings of $\mathbb{P}^1$ into $\mathbb{P}^4$ that map birationally
onto images of degree $10$, of arithmetic genus $g$, and of span an
$r$-dimensional projective space. Note that we do {\it not} mod out
by the $\mbox{PGL}(2)$-action on $\mathbb{P}^1$; nor do we projectivize. Let $I_{10}^{r,g}$ denote the incidence locus in
$M_{10}^{r,g} \times \mathbb{P}^{125}$. Let $M_{10,i}^{r,g} \subset
M_{10}^{r,g}$ denote the sublocus of morphisms with images $C$
satisfying
\begin{equation*}
h^1(\mathcal{I}_C(5))=i, 
\end{equation*}
and let $I_{10,i}^{r,g} \subset I_{10}^{r,g}$ denote
the pullback of $M_{10,i}^{r,g}$ under the canonical projection $I_{10}^{r,g}
\rightarrow \mathbb{P}^{125}$.

As explained in \cite{Klei}, to
prove Theorem~\ref{1.1} it's enough to show that the projection $I_{10}^{r,g}
\rightarrow \mathbb{P}^{125}$ is surjective only over the irreducible
component $M_{10,0}$ of morphisms whose images $C$ are smooth with
$h^1(\mathcal{I}_{C}(5))=0$. In particular, this
statement implies
the finiteness of the Hilbert scheme of curves on a general quintic,
and the splitting property of the normal bundles
$\mathcal{N}_{C/F}$ follows immediately from Verdier's result \cite[Thm,
  p.139]{Ve}. See the proof of \cite[Cor
  2.5]{Klei} for more details.

However, since the projection's fibres all have dimension at least 4,
to prove that it's not surjective over the locus of curves with
$i$ nonzero, it
suffices to establish that $\dim I_{10}^{r,g}<129$. Using
the exact
sequence
\begin{equation*}
0 \rightarrow H^0(\mathcal{I}_C(5)) \rightarrow
H^0(\mathcal{O}_{\mathbb{P}^4}(5)) \rightarrow H^0(\mathcal{O}_C(5))
\rightarrow H^1(\mathcal{I}_C(5)) \rightarrow 0,
\end{equation*}
to relate the Hilbert polynomial of $\mathcal{I_C}$ to that of $C$, we see that
\begin{equation*}
\dim I_{10,i}^{r,g} \leq \dim M_{10,i}^{r,g}+126-(5 \cdot 10-g+1)+i,
\end{equation*}
at least provided $H^1(\mathcal{O}_C(5))=0$. (See \cite{Klei} for more
details.) Assuming this vanishing, which we prove in Lemma~\ref{L2} below, we are reduced to showing that
\begin{equation}\label{55-g-i}
\dim M_{10,i}^{r,g} < 55-g-i
\end{equation}
whenever $i$ is nonzero.

In order to simplify the notation, we let $i:=h^1(\mathcal{I}_C(5))$
hereafter.

We may also assume $r \geq 3$, since (as noted in \cite[Section 3]{Klei})
no planar rational curves of degree greater than 5 lie on a general quintic.

Moreover, by the result of \cite[Lemma 3.4]{Klei}, we have
\begin{subequations}\label{1.2}
\begin{align}
\dim M_d^{4,g} &\leq 5(d+1)-1-\min(2g,8)\label{1.2a}\text{ for all }d
\geq 7, \text{ and }\\
\dim M_d^{3,g} &\leq 4(d+1)-1-\min(g,5)\label{1.2b}\text{ for all }d
\geq 9.
\end{align}
\end{subequations}

To establish the validity of \eqref{55-g-i}, we treat nondegenerate
curves and curves with 3-planar images separately. It's interesting to
note that in this paper, unlike in \cite{Klei}, the main technical work lies in the proof of
\eqref{55-g-i} in the
nondegenerate case. Accordingly, we first consider nondegenerate
curves in $\mathbb{P}^4$.
Applying \eqref{1.2a}, it's easy to see that Theorem~\ref{1.1} for nondegenerate
curves amounts to the following.

\begin{thm}\label{1.1bis}
The nondegenerate reduced irreducible
rational curves verifying
\begin{equation*}
g+i \geq \min(2g,8)
\end{equation*}
define a sublocus of the parameter space $M_{10}$ of degree-$10$
rational maps of codimension greater than
$g+i$.
\end{thm}

\subsection{Special subloci of the mapping space}
In this subsection, we obtain lower bounds on the codimensions of
certain special subloci of $M_{10}$. These estimates effectively
allow us to ignore certain curves $C$ corresponding to generic initial ideals with
large $h^1(\mathcal{I}_C(5))$.

\begin{Lemma}\label{1a}
Degree-$10$ morphisms determining
nondegenerate degree-$10$ rational curves
with $8$-secant lines have codimension at least $10$ in the space of rational
mappings of degree $10$. 
\end{Lemma}

\begin{proof} Every rational map $f\:\mathbb{P}^1
\rightarrow \mathbb{P}^4$ of degree $10$ is parametrized by 5 equations
$x_i=f_i(t,u) \text{ for }i=0,\dots,4$, where $(t,u)$ are coordinates on
$\mathbb{P}^1$, the $x_i$ are coordinates on $\mathbb{P}^4$, and $f_i$
is a polynomial of degree $10$ in
$t$ and $u$. To say that the image of $f$ intersects the line
$L=\{x_0=x_1=x_2=0\}$ in {\it distinct}
points $f(t_j,u_j)$ for $1 \leq j \leq 8$,  means precisely that
\begin{equation}\label{f_i}
f_i(t_j,u_j)=0 \text{ for } 0 \leq i \leq 2 \text{ and } 1
\leq j \leq 8.
\end{equation}  
In other words, $24$ linearly independent equations must be satisfied. So requiring a
map $f$ of degree $10$ to map any fixed degree-$8$ divisor onto distinct
points along any fixed
line $L$ is a codimension-$24$ condition. Moreover, if we allow the
image points to coalesce and replace the equations \eqref{f_i} by
their higher-contact analogues, then the same argument shows that requiring $f$ of degree $10$ to map any fixed degree-$8$ divisor along any fixed
line $L$ is a codimension-$24$ condition. On the other hand, the degree-$8$ divisors on $\mathbb{P}^1$
determine a $\mathbb{P}^8$ and the Grassmannian of lines in
$\mathbb{P}^4$ has dimension $6$, so by varying our choice of eight
points in $\mathbb{P}^1$ along with the choice of the line $L \subset
\mathbb{P}^4$, we obtain $10$ conditions and the result
follows.
\end{proof}

\begin{Lemma}\label{1b}Degree-$10$ morphisms defining rational curves in reduced, irreducible hyperquadrics determine a locus of codimension
  at least 7 in 
$M_{10}$.
\end{Lemma}

\begin{proof} By the main result of \cite{KiPa}, the
  space of rational curves of any fixed degree $d$ inside a projective
  homogeneous space is irreducible and of the expected dimension. On the other hand, every {\it smooth} hyperquadric is a
homogeneous space, and any two hyperquadrics are isomorphic, so the codimension of curves lying on
smooth hyperquadrics is the expected one. Note that hyperquadrics $F \subset \mathbb{P}^4$
comprise a $\mathbb{P}^{14}$ and that the equation defining the
incidence scheme
$\{C \subset F: C \in M_{10}, F \in \mathbb{P}^{14}$ in
$M_{10} \times \mathbb{P}^{14}\}$ is a polynomial in the coordinates
of $\mathbb{P}^1$ in $21$ parameters. So $21$ linearly independent
conditions must be met for
a degree-$10$ curve to lie on a fixed smooth hyperquadric. By varying the
choice of hyperquadric, we find that the codimension of (morphisms
corresponding to) curves inside smooth hyperquadrics is $7$.

Now let $Q$ be a singular, reduced, and irreducible hyperquadric. Since it contains a
 nondegenerate curve $C$, it must be a cone over a
 quadric $Q$ of rank $2$ or $3$. 

Say $\mbox{rank}(Q)=3$. Let
$\mbox{Hom}_{10}(\mathbb{P}^1,Y)$ denote the (affine) parameter space of morphisms $\mathbb{P}^1 \rightarrow Y$. When $Y=\mathbb{P}^4$, a morphism $f:\mathbb{P}^1
\rightarrow Y$ is given by $5$ sections
$f_i \in H^0(\mathcal{O}_{\mathbb{P}^1}(\deg(f)))$ for $i=0,\dots,4$. Projection from the
vertex of $Q$ defines a rational map
$\pi:\mbox{Hom}_{10}(\mathbb{P}^1,\mathbb{P}^4) \rightarrow
\mbox{Hom}_{10}(\mathbb{P}^1,\mathbb{P}^3)$. Now assume that the vertex $p$ of $Q$ has coordinates
$(0,0,0,0,1)$. Then $\pi: \mbox{Hom}_{10}(\mathbb{P}^1,\mathbb{P}^4)$
is the rational map that drops the fifth parameterizing polynomial:
\begin{equation*}
(f_0,f_1,f_2,f_3,f_4) \mapsto (f_0,f_1,f_2,f_3).
\end{equation*}
Since $f_4$ has $11$ coefficients, the (affine)
fibre dimension of $\pi$ is $11$ over any degree-$10$ map $f$ whose image curve avoids
$p$. Any such map $f \in
\mbox{Hom}_{10}(\mathbb{P}^1,Q)$ is sent under
$\pi$ to a map $\stackrel{\sim}{f} \in
\mbox{Hom}_{10}(\mathbb{P}^1,\stackrel{\sim}{Q})$, where
$\stackrel{\sim}{Q} \cong \mathbb{P}^1 \times \mathbb{P}^1$ is a smooth
quadric surface. As explained in \cite[Lemma 3.2]{Klei}, the scheme associated to $\mbox{Hom}_{10}(\mathbb{P}^1,\stackrel{\sim}{Q})$ is a disjoint
union of
$23$-dimensional projective spaces.

On the other hand, quadric cones determine a divisor inside the
$\mathbb{P}^{14}$ of hyperquadrics in $\mathbb{P}^4$. Since rational
maps whose images pass through the vertex of $Q$ clearly comprise a
proper closed subset of $\mbox{Hom}_{10}(\mathbb{P}^1,Q)$, we conclude that rational maps to cones over quadric
surfaces vary in a family of
dimension at most $13+11+23=47$. In other words, such maps determine a locus of
codimension $55-47=8$ inside $M_{10}^4$.

To complete the proof of Lemma~\ref{1b}, we must handle rational curves
whose images lie in cones $Q^{\prime}$ with 1-dimensional vertices over plane
conics. In fact, these pose no problem. For, by a result of
\cite{Gu} quoted in \cite[proof of Lemma 3.2]{Klei},
$\mbox{Hom}_{10}(\mathbb{P}^1,Q^{\prime})$ is at most
$23$-dimensional. We now conclude by applying the same
projection-based argument used in the preceding case.
\end{proof}

\subsection{A useful stratification of the mapping space}
Let $f\:\mathbb{P}^1 \rightarrow \mathbb{P}^4$ be a map with image
$C$. Note that the $(-1)$-twist of the ``restricted
tangent bundle'' $f^* T_{\mathbb{P}^4}(-1)$, has degree $10$ and rank $4$ over
$\mathbb{P}^1$, say with splitting 
\begin{equation}
f^* T_{\mathbb{P}^4}(-1)=
\bigoplus_{i=1}^4 \mathcal{O}_{\mathbb{P}^1}(a_i). \notag
\end{equation}
A result of Verdier's (see \cite{Ve} and \cite{Ra})
establishes that the scheme of morphisms $\mathbb{P}^1 \rightarrow
\mathbb{P}^4$ of fixed degree corresponding to a particular splitting type
$(a_1,a_2,a_3,a_4)$ with $a_1 \geq a_2 \geq a_3 \geq a_4$ is
irreducible of the expected codimension
\begin{equation}
\sum_{i\neq j} \max\{0,a_i-a_j-1\}. \notag
\end{equation}

In particular, every special splitting stratum has
codimension $4$ or more, and the stratum $(4,3,2,1)$ is the unique
stratum of codimension $4$. If $f^* T_{\mathbb{P}^4}(-1)$ has
the generic splitting, then $\bigwedge^2 (f^* T_{\mathbb{P}^4}(-1))^{*} \otimes
\mathcal{O}_{\mathbb{P}^1}(5)$ has no higher cohomology, and it follows by
\cite[Proposition 1.2]{GLP} that $\mathcal{I}_C$ is $6$-regular. So $i \neq 0$ in codimension
$4$. Therefore, in light of the discussion immediately preceding the
statement of Theorem~\ref{1.1bis}, we have reduced to showing that for all
nondegenerate $C$ such that $i \neq 0$, one of following conditions is verified:
\begin{subequations}
\begin{align}
&g+i<4, \label{1.3a} \text{ or }\\
&C \text{ admits an 8-secant line and }g+i<10, \label{1.3b}\text{ or }\\
&C \text{ lies on a hyperquadric and }g+i<7, \label{1.3c}\text{ or }\\
&g+i<\min(2g,8). \label{1.3d}
\end{align}
\end{subequations}
A curve $C$ meeting any of these conditions will be called
{\it nonproblematic}. Similarly, any sublocus of $M_{10}$ comprised of
nonproblematic curves will be called nonproblematic.

\subsection{Proof of Theorem~\ref{1.1} for nondegenerate curves}
To prove Theorem~\ref{1.1}, we obtain bounds on
$h^1(\mathcal{I}_C(5))$, based upon an analysis of possible
corresponding Borel-fixed monomial ideals. Each
$\mbox{gin}(\mathcal{I}_C)$ is obtained, by Fact $1$, by applying a
sequence of rewriting rules to $\mbox{gin}(\mathcal{I}_{\Gamma})$. We
now argue case by case, based upon the possibilities for $\mbox{gin}(\mathcal{I}_{\Gamma})$.

\newtheorem{Case}{Case}
\begin{Case}\label{C1}
Say $\mbox{gin}(\mathcal{I}_{\Gamma})= \mbox{Borel}(x_2^3)$.
\end{Case}
Because
$I_{C_{\Gamma}}$ is minimally generated by polynomials of degree at
most $3$, $C_{\Gamma}$ is $3$-regular. So Lemma \ref{g_Gamma} implies that
\begin{equation*}
g_{\Gamma}=h^0(\mathcal{I}_{C_{\Gamma}}(3))-4. 
\end{equation*}
But, because $I_{C_{\Gamma}}$ is saturated, $h^0(\mathcal{I}_{C_{\Gamma}}(3))$ is equal to the number
of linearly independent polynomials of degree $3$ in
$I_{C_{\Gamma}}$. Since there are exactly ten of these, we find $g_{\Gamma}=6$.

It follows from Lemma \ref{g+i} that
$g+i \leq 6$. Note that if $C$
as above admits an $8$-secant line, then $C$ belongs to a locally-closed
subscheme inside $M_{10}$ of codimension at least $10$, by Lemma
\ref{1a}. Thus $C$ is
nonproblematic. The remaining curves are $7$-regular, by
\cite[Thm.~3.1, p.501]{GLP}; we will see that
they, too, are nonproblematic. For this purpose, we begin by showing
that each one satisfies 
\begin{equation}\label{iboundC1}
i \leq 1.
\end{equation}
By Lemma \ref{g+i}, $\mbox{gin}(\mathcal{I}_C)$ may be obtained
  from $I_{C_{\Gamma}}$ by
applying
precisely $6-g$ rewriting rules. On the other hand, by Lemma~\ref{4}, $i$ is exactly the number of
$C$-rewritings applied in degrees $6$ or greater. In particular, $i$
  is maximized when the number of $C$-rewritings (of arbitrary
  degrees) is maximized, which happens when $g=0$. Therefore, for the
  purpose of verifying \eqref{iboundC1}, we assume $g=0$. 

There are now a large number of possible $C$-rewriting sequences to choose
from. For the sake of bounding $i$, however, we only
need consider those rewriting sequences that result in the maximal number
of minimal generators of degree greater than six. Moreover, certain
rewritings are dictated by Borel-fixity. If $i$ is nonzero, for
instance, then the minimal generator $x_2^3$ {\it must} be
rewritten. Without loss of generality, therefore, we assume the first rewriting rule exchanges the
minimal generator $x_2^3$ for minimal generators $x_2^4$ and
$x_2^3x_3$, i.e., that it is $x_2 \mapsto (x_2,x_3)$
applied at the leaf corresponding to the minimal generator $x_2^3$ of
$\mbox{gin }\mathcal{I}_{\Gamma}$. 

\begin{figure}
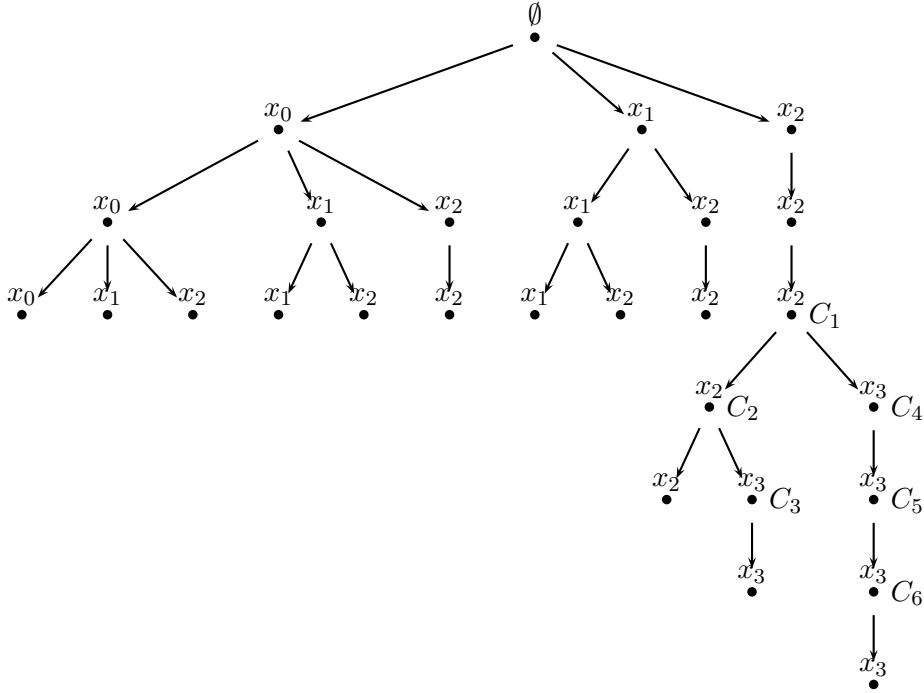

\begin{center}
\psset{arrows=->,labelsep=1pt,nodesep=7pt,tnpos=a,levelsep=35pt}

\pstree{\Tdot*~{$\emptyset$}}{

\pstree{\Tdot*~{$x_0$}}{ 

\pstree{\Tdot*~{$x_0$}} 
{\Tdot*~{$x_0$} \Tdot*~{$x_1$} \Tdot*~{$x_2$}}

\pstree{\Tdot*~{$x_1$}} 
{\Tdot*~{$x_1$} \Tdot*~{$x_2$}}

\pstree{\Tdot*~{$x_2$}} 
{\Tdot*~{$x_2$}}
}

\pstree{\Tdot*~{$x_1$}}{ 

\pstree{\Tdot*~{$x_1$}} 
{\Tdot*~{$x_1$} \Tdot*~{$x_2$}}

\pstree{\Tdot*~{$x_2$}} 
{\Tdot*~{$x_2$}}
}

\pstree{\Tdot*~{$x_2$}}{ 

\pstree{\Tdot*~{$x_2$}}{ 

\pstree{\Tdot*~{$x_2$}~[tnpos=r]{$C_1$}}{ 

\pstree{\Tdot*~{$x_2$}~[tnpos=r]{$C_2$}}{ 
\Tdot*~{$x_2$} \pstree{\Tdot*~{$x_3$}~[tnpos=r]{$C_3$}}{
\Tdot*~{$x_3$}}}
\pstree{\Tdot*~{$x_3$}~[tnpos=r]{$C_4$}}{ 
\pstree{\Tdot*~{$x_3$}~[tnpos=r]{$C_5$}}{
\pstree{\Tdot*~{$x_3$}~[tnpos=r]{$C_6$}}{
\Tdot*~{$x_3$}}
}}}

}

}

}
\caption{$\mbox{gin
  }\mathcal{I}_C=\mbox{Borel}(x_1x_2^2)+(x_2^5,x_2^4x_3^2,x_2^3x_3^4)$.
  The $C$-rules are marked as $C_i, i=1,\dots,6$.}
\label{figure: Figure $6$}

\end{center}
\end{figure}

Similarly, we may assume the next
three rewriting rules have the effect of exchanging the minimal
generator $x_2^3x_3$ for $x_2^3x_3^4$. But then no further rewriting
may be applied to $x_2^3x_3^4$, since $C$ is $7$-regular. So
rewriting twice more, beginning with minimal generators of degree at
most $4$, results in minimal generators of degree at most $6$; for
instance, we might have
\begin{equation}
\mbox{gin}(\mathcal{I}_C)=\mbox{Borel}(x_1x_2^2)+(x_2^5,x_2^4x_3^2,x_2^3x_3^4). \notag
\end{equation}
The latter ideal is obtained from the (extension of the) hyperplane gin by applying
$C$-rules to leaves that are ``farthest to the right'' in the tree of Figure 2;
see Figure~\ref{figure: Figure $6$}. It follows immediately by Lemma~\ref{4} that $i \leq 1$.

It follows that
\begin{equation*}
g+i < 4
\end{equation*}
whenever $g \leq 2$. Therefore, if $g \leq 2$, then $C$ is
nonproblematic.

It remains to show that $7$-regular $C$ with $g \geq 3$ and $g+i \leq
6$ are nonproblematic. By Lemma~\ref{g+i}, every such $C$ is such
that $\mbox{gin}(\mathcal{I}_C)$ is obtained from the hyperplane gin
by applying at most $3$-rewriting rules, all of which are necessarily
in degrees less than six. By Lemma~\ref{4} it follows that $i=0$. So every such $C$ is indeed nonproblematic.

\begin{Case}\label{C2}
Say $\mbox{gin
}\mathcal{I}_{\Gamma}=\mbox{Borel}(x_1x_2^2)+(x_0^2,x_2^4)$.
\end{Case} 
Because
$I_{C_{\Gamma}}$ is minimally generated by polynomials of degree at
most $4$, $C_{\Gamma}$ is $4$-regular. So Lemma \ref{g_Gamma} implies that
\begin{equation*}
g_{\Gamma}=h^0(\mathcal{I}_{C_{\Gamma}}(4))-29. 
\end{equation*}
Because $I_{C_{\Gamma}}$ is saturated, $h^0(\mathcal{I}_{C_{\Gamma}}(4))$ is equal to the number
of linearly independent polynomials of degree $4$ in
$I_{C_{\Gamma}}$; there are $36$ of these, so $g_{\Gamma}=7$.
It follows from Lemma \ref{g+i} that
$g+i \leq 7$. It follows, as in Case \ref{C1}, that if $C$
admits an $8$-secant line, then $C$ is
nonproblematic. The remaining curves are $7$-regular, by
\cite[Thm.~3.1, p.501]{GLP}; we now show they are also
nonproblematic.

By Lemma~\ref{4}, $i$ is exactly the number of
$C$-rewritings applied in degrees $6$ or greater. But by Lemma \ref{g+i}, $\mbox{gin}(\mathcal{I}_C)$ may be obtained
  from $I_{C_{\Gamma}}$ by
applying
precisely $7-g$ rewriting rules. But because $C$ is $7$-regular, it
follows there are no rewritings in degree greater than $6$. Moreover,
there are at
most three rewritings in degree $6$, because it takes six
rewritings total to obtain three rewritings in degree $6$, and any
seventh rewriting will be in degree less than $6$. (Compare Figure $3$.) Therefore, $i \leq 3$.

It follows immediately that if $g=0$, then $g+i \leq 3$ and,
therefore, $C$ is nonproblematic.

\begin{figure}
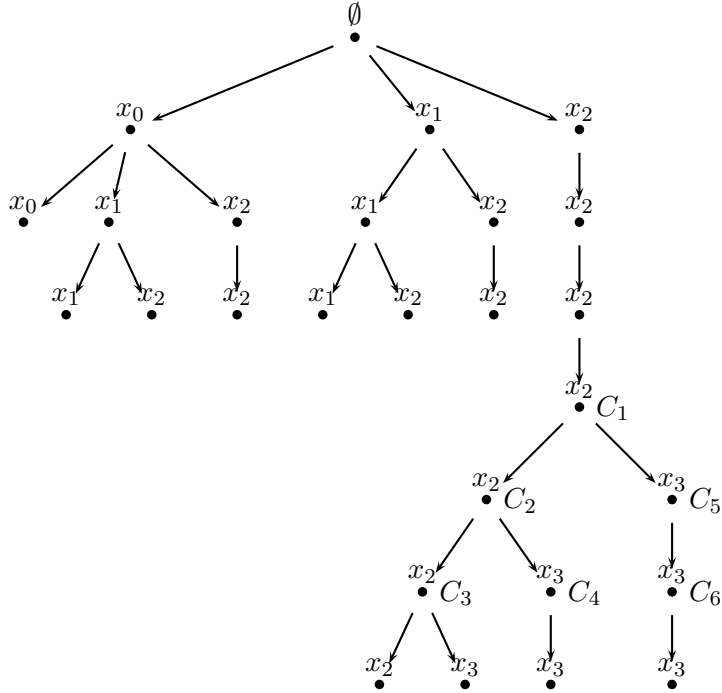

\begin{center}
\psset{arrows=->,labelsep=1pt,nodesep=7pt,tnpos=a,levelsep=35pt}

\pstree{\Tdot*~{$\emptyset$}}{

\pstree{\Tdot*~{$x_0$}}{ 

\Tdot*~{$x_0$}

\pstree{\Tdot*~{$x_1$}} 
{\Tdot*~{$x_1$} \Tdot*~{$x_2$}}

\pstree{\Tdot*~{$x_2$}} 
{\Tdot*~{$x_2$}}
}

\pstree{\Tdot*~{$x_1$}}{ 

\pstree{\Tdot*~{$x_1$}} 
{\Tdot*~{$x_1$} \Tdot*~{$x_2$}}

\pstree{\Tdot*~{$x_2$}} 
{\Tdot*~{$x_2$}}
}

\pstree{\Tdot*~{$x_2$}}{ 

\pstree{\Tdot*~{$x_2$}}{ 

\pstree{\Tdot*~{$x_2$}}{ 
\pstree{\Tdot*~{$x_2$}~[tnpos=r]{$C_1$}}{ 

\pstree{\Tdot*~{$x_2$}~[tnpos=r]{$C_2$}}{ 
\pstree{\Tdot*~{$x_2$}~[tnpos=r]{$C_3$}}{\Tdot*~{$x_2$}
  \Tdot*~{$x_3$}} 
\pstree{\Tdot*~{$x_3$}~[tnpos=r]{$C_4$}}{
\Tdot*~{$x_3$}}}
\pstree{\Tdot*~{$x_3$}~[tnpos=r]{$C_5$}}{ 
\pstree{\Tdot*~{$x_3$}~[tnpos=r]{$C_6$}}{
\Tdot*~{$x_3$}}
}}

}

}

}
}
\caption{$\mbox{gin
  }\mathcal{I}_C=(x_0^2)+\mbox{Borel}(x_1x_2^2)+(x_2^7,x_2^6x_3,x_2^5x_3^2,x_2^4x_3^3)$.
  The $C$-rules are marked as $C_i, i=1,\dots,6$. Clearly,
  $g=1$. Moreover, $i=3$, since there are
  exactly three $C$-rules in degree 6.}
\label{figure: Figure $7$}

\end{center}
\end{figure}

Similarly, if $g=1$, then either $g+i \leq 3$, in which case $C$  is
nonproblematic, or $g+i=4$ and $\mbox{gin}(\mathcal{I}_C)$ corresponds
to the tree of Figure~\ref{figure: Figure $7$}. In the latter situation, however, Figure
$7$ shows that $C$ lies
on a hyperquadric, and therefore, by Lemma~\ref{1b}, $C$ is nonproblematic.

Similarly, if $g=2$, then $\mbox{gin}(\mathcal{I}_C)$ is obtained from the
hyperplane gin in $5$ rewritings, of which at most $2$ may be in
degree $6$. So $i \leq 2$ by Lemma~\ref{4}, i.e., $g+i \leq 4$. Moreover, equality is
obtained if and only if $C$ lies on a hyperquadric, in which case $C$
is nonproblematic, by Lemma~\ref{1b}.

Finally, if $g \geq 3$, then $\mbox{gin}(\mathcal{I}_C)$ is obtained from the
hyperplane gin in at most $4$ rewritings, so $i \leq 2$, which implies
$g+i<\min(2g,8)$. So $C$ is nonproblematic in this situation as well,
which enables us to conclude.

\begin{Case}\label{C3}
Say $\mbox{gin }\mathcal{I}_{\Gamma}=\mbox{Borel
}(x_0x_1)+(x_1^3,x_1^2x_2,x_0x_2^2)+(x_1x_2^3,x_2^4)$.
\end{Case}

To show that the corresponding curves $C$ are nonproblematic,
we proceed as follows.

First, $C$ is $4$-regular, so 
\begin{equation*}
g_{\Gamma}=h^0(\mathcal{I}_{C_{\Gamma}}(4))-29,
\end{equation*}
by Lemma \ref{g_Gamma}. Since $I_{C_{\Gamma}}$ contains $37$ linearly
independent polynomials in degree $4$, we have $g_{\Gamma}=8$. Hence, $g+i
\leq 8$, by Lemma \ref{g+i}. Note that applying $r \leq 2$ $C$-rewriting
rules to $I_{C_{\Gamma}}$ results in a Borel-fixed ideal
that is
minimally generated in degrees at most $6$ and has genus
$g_{\Gamma}-r$, by Lemma \ref{L3}. So if $g=6$ or $g=7$, then
$\mbox{gin }\mathcal{I}_C$ is $6$-regular, by \cite[Thm.~2.27]{Green},
and moreover
\begin{equation*}
g+i < \min(2g,8),
\end{equation*}
and therefore $C$ is $6$-regular. We now show curves of genus 8
are nonproblematic, too. 

Note that any such curve $C$ with $g=8$
necessarily satisfies
\begin{equation*}
\mbox{gin}(\mathcal{I}_C)=I_{C_{\Gamma}},
\end{equation*}
by Lemma~\ref{g+i}. Since $h^0(\mathcal{I}_{C_{\Gamma}}(2))=2$, it follows that $C$ lies on two linearly independent hyperquadrics
$Q_1$ and $Q_2$. These intersect in a
surface $S$ of degree $4$. In fact, $C$ is also contained in a
hypercubic, $K$, that is linearly independent of $Q_1$ and $Q_2$. To
see why, assume that no such $K$ exists. Then every hypercubic $K$
containing $C$ satisfies $Q_1 \cap Q_2 \subset K$, so we must have
\begin{equation*}
h^0(\mathcal{I}_C(3)) \leq h^0(\mathcal{I}_{Q_1 \cap Q_2}(3)),
\end{equation*}
and therefore
\begin{equation}\label{sections_ineq}
h^0(\mathcal{I}_{\Gamma}(3)) \leq h^0(\mathcal{I}_{\Lambda}(3)),
\end{equation}
where $\Gamma$ and $\Lambda$ denote general hyperplane sections of $C$
and $Q_1 \cap Q_2$, respectively.
Note that $Q_1$ and $Q_2$ are nondegenerate, since they contain $C$;
so $Q_1 \cap Q_2$ is a nondegenerate quartic surface, and $\Lambda$ is
a nondegenerate quartic space curve. By \cite[Thm.~3.1, p.501]{GLP},
$\Lambda$ is $3$-regular, and it follows from the usual long exact
sequence in cohomology that $\mathcal{O}_{\Lambda}$ is $3$-regular,
too. So there is an exact sequence
\begin{equation*}
0 \rightarrow H^0(\mathcal{I}_{\Lambda}(3)) \rightarrow
H^0(\mathcal{O}_{\mathbb{P}^3}(3)) \rightarrow
H^0(\mathcal{O}_{\Lambda}(3)) \rightarrow 0,
\end{equation*}
and by the Riemann-Roch theorem, we deduce that
\begin{equation*}
h^0(\mathcal{O}_{\Lambda}(3)) = 13-g(\Lambda).
\end{equation*}
By \cite{Ciliberto}, $\Lambda$ has genus at most 1; it follows that
\begin{equation*}
h^0(\mathcal{I}_{\Lambda}(3)) \leq 8.
\end{equation*}
On the other hand, we have
\begin{equation*}
h^0(\mathcal{I}_{\Gamma}(3))=10,
\end{equation*}
so \eqref{sections_ineq} is violated, which gives a contradiction.

Now let $S:=Q_1 \cap Q_2$. If $S \cap K$ contains a surface component $S^{\star}$,
then clearly $S^{\star}$ is properly contained in $S$, whence
$\deg(S^{\star})<4$. Note that $C$ lies on either $S^{\star}$ or on a
component $S^{\star \star}$ belonging to the residual to $S^{\star}$ in $Q_1 \cap
Q_2$. However, because $C$ is nondegenerate, $C$ lies on no component
of $S$ degree $1$ or $2$. So the component of $S$ on which $C$ lies is a nondegenerate threefold of degree $3$,
which is necessarily a cubic scroll by \cite[Prop.,
  p.525]{GH}. Moreover, no cubic scroll containing
$C$ meets $K$ properly, since $\deg(C)=10$. So $S^{\star \star}$ is a
surface component of $S \cap K$. Therefore, replacing $S^{\star}$ by $S^{\star
  \star}$ if necessary, we may assume $S^{\star}$ is a cubic
scroll containing $C$.

We handle degree-$10$ curves on cubic scrolls in $\mathbb{P}^4$
as follows. As explained in \cite[pp.519-523]{GH}, cubic scrolls come
in two basic types. Another good reference for scrolls, whose
notation we will use and to
which we will often refer, is \cite{Cos}. Those that are
smooth will be denoted $S_{1,2}$ (in the notation of \cite{GH}, these
are the scrolls $S_{1,1}$). Those that are singular are
cones over twisted cubic curves, and will be denoted by $S_{0,3}$. 

Every $S_{1,2}$ may be realized as the
image of
the Hirzebruch surface $F_1:=\mathbb{P}(\mathcal{O}_{\mathbb{P}^1} \oplus
\mathcal{O}_{\mathbb{P}^1}(1))$ under the map $\phi_{1,2}$ defined by the complete linear series
$|O_{F_1}(e+2f)|$, where $e$ is the divisor class of the
$-1$-curve on $F_1 \rightarrow \mathbb{P}^1$ and $f$ is the class of
the fibre.

The intersection pairing on $F_1$ is given by
\begin{equation*}
e^2=-1, e \cdot f=1, f^2=0.
\end{equation*}
Moreover, the
canonical class of $F_1$ is
\begin{equation*}
K_{F_1}=-2e-3f.
\end{equation*}

Let $\phi_{1,2}^{\star}[C]=ae+bf \in \mbox{Pic }F_1$. The adjunction formula implies
that the arithmetic genus of $C$ satisfies
\begin{equation}\label{adjunction_on_S_{1,2}}
2g-2=((a-2)e+(b-3)f)\cdot(ae+bf),
\end{equation}
i.e.,
\begin{equation}\label{S_{1,2}eq1}
a^2-2ab+a+2b+2g-2=0.
\end{equation}
Also, since $\deg(C)=10$, we have
\begin{equation*}
(e+2f) \cdot (ae+bf)=10,
\end{equation*}
or equivalently,
\begin{equation}\label{S_{1,2}eq2}
b=10-a.
\end{equation}
Substituting \eqref{S_{1,2}eq2} in \eqref{S_{1,2}eq1}, we obtain
\begin{equation}\label{S_{1,2}eq3}
3a^2-21a+18+2g=0.
\end{equation}
If $g=8$, then \eqref{S_{1,2}eq3} has no integral solutions $a$.

Finally, we treat curves $C$ lying on cones $S_{0,3}$ over twisted
cubics. Every $S_{0,3}$ is the image of the Hirzebruch
surface $F_3=\mathbb{P}(\mathcal{O}_{\mathbb{P}^1} \oplus \mathcal{O}_{\mathbb{P}^1}(3))$
under the map $\phi_{0,3}$ defined by the complete linear series $|O_{F_3}(e+3f)|$, where $e$ is the
divisor class of the  of $-3$-curve on $F_3$ and $f$ is the class of the fibre. Note that $\phi_{0,3}$ is a
birational map that blows down the $-3$-curve to the vertex
of $S_{0,3}$.

The intersection pairing on $F_3$ is given by
\begin{equation*}
e^2=-3, e \cdot f=1, f^2=0.
\end{equation*}
The
canonical class of $F_3$ is
\begin{equation}\label{K_F3}
K_{F_3}=-2e-5f.
\end{equation}

Now say that $C$
avoids the vertex $p$ of $S_{0,3}$. Then
$[\phi_{0,3}^{-1}C]=ae+bf$, for some integers $a$ and $b$. In fact,
since $\deg(C)=10$, we have
\begin{equation*}
(e+3f) \cdot (ae+bf)=10,
\end{equation*}
whence,
\begin{equation}
b=10.
\end{equation}

On the other hand, the adjunction formula
then implies that
\begin{equation}\label{S_{0,3}eq2}
2g-2=((a-2)e+(b-5)f)\cdot (ae+bf).
\end{equation}
Substituting $b=10$ in \eqref{S_{0,3}eq2} and solving for $a$ yields
\begin{equation}\label{S_{0,3}eq1}
a=\frac{2}{7}g+1.
\end{equation}
Here $g=8$; by \eqref{S_{0,3}eq1}, it follows $a$ is not an integer,
which is plainly absurd. 

Now say that $C$ passes through the vertex of
$S_{0,3}$ with multiplicity $m>0$.
Then $[\phi_{0,3}^{-1}C]=ae+bf$, for some integer $b$, and
\begin{equation*}
[\stackrel{\sim}{C}]=(a-m)e+bf,
\end{equation*}
where
$\stackrel{\sim}{C}$ denotes the proper transform of $C$ on
$F_3$. 

Let
$\stackrel{\sim}{g}:=g(\stackrel{\sim}{C})$. Then we have
\begin{equation}\label{eq1}
2\stackrel{\sim}{g}-2=((a-m-2)e+(b-5)f)\cdot ((a-m)e+bf)
\end{equation}
by the adjunction formula,
\begin{equation}\label{eq2}
\stackrel{\sim}{g} \leq g,
\end{equation}
and
\begin{equation*}
(e+3f) \cdot (ae+bf)=10,
\end{equation*}
since $\deg(C)=10$.
The latter equation implies
\begin{equation*}
b=10,
\end{equation*}
just as before.

Substituting $b=10$ in
\eqref{eq1} and solving for $a-m$ yields 
\begin{equation}\label{eq4}
a-m=\frac{2}{7}\stackrel{\sim}{g}+1,
\end{equation}
for some $\stackrel{\sim}{g} \leq 8$. Because $a$ and $m$ are integers
it follows from \eqref{eq4} that 
\begin{equation*}
\stackrel{\sim}{g}=7 \text{ and }a-m=3.
\end{equation*}
In particular, the genus discrepancy $g-\stackrel{\sim}{g}$ is equal to $1$; whence,
\begin{equation*}
m \leq 2,
\end{equation*}
by the result of Rosenlicht quoted in \cite[Ch. VIII, Prop. 1.16]{AK}.
Therefore, there are two possibilities for the class of $[\phi_{0,3}^{-1}C]$:
either
\begin{equation*}
[\phi_{0,3}^{-1}C]=4e+10f,
\end{equation*}
or
\begin{equation*}
[\phi_{0,3}^{-1}C]=5e+10f.
\end{equation*}

To handle curves of the latter two types, we argue as follows. Since $\mathbb{F}_{3}$
is rational, we have $\chi(O_{\mathbb{F}_{3}})=1$. Therefore, the Riemann-Roch 
theorem for surfaces implies that
\begin{equation}
\begin{split}
\chi(ae+bf)&=1+\frac{1}{2}((ae+bf)^{2}-(ae+bf)\cdot(-2e-5f)). \\
&=-\frac{3}{2}a^{2}+ab-a+ 2b+1.
\end{split} \notag
\end{equation}
Substituting $(a,b)=(4,10)$ and $(a,b)=(5,10)$ yields
\begin{equation*}
\chi(ae+bf)=33       
\end{equation*}
and
\begin{equation*}
\chi(ae+bf)=\frac{57}{2},
\end{equation*}
respectively. However, $\chi(ae+bf)$ is necessarily an integer, so we must have 
\begin{equation*}
(a,b)=(4,10).
\end{equation*}
Then \eqref{eq4} implies that the proper transform of $C$ on $\mathbb{F}_{3}$
has class
\begin{equation*}
[\stackrel{\sim}{C}]=e+10f.
\end{equation*}
However, the Riemann-Roch formula then implies that $\chi(\stackrel{\sim}{C})$
is not an integer, which is absurd.

Now consider the only remaining possibility, that $S \cap K$ is a complete
intersection $X$ of type $(2,2,3)$. Thus $C$ is residual in $X$ to a
$1$-dimensional scheme $R$ of degree $2$. Let $\mathcal{L}$ denote the
linear series of hypercubics containing $C$. By Bertini's theorem, a
general member of $\mathcal{L}$ is smooth away from $\mbox{Sing}(C)
\cup \mbox{Bs}(\mathcal{L})$, where $\mbox{Bs}(\mathcal{L})$ denotes
the base locus of $\mathcal{L}$. Assuming, as we may, that the
quadrics and cubics defining $X$ are general among quadrics and cubics
containing $C$, it follows that $R$ is reduced (i.e., a plane conic), since otherwise
\begin{equation*}
X \subset \mbox{Bs}(\mathcal{L}),
\end{equation*}
which is absurd. On the other hand, the adjunction formula implies that
\begin{equation*}
K_X=\mathcal{O}_X(2) \text{ and }g(X)=13,
\end{equation*}
So $g(R)=16$,
by the main result of \cite{Na}, which is also absurd.

It follows that every curve $C$ of genus $8$ is
  nonproblematic. 

Next, say $g=6 \text{ or 7}$. Then $\mbox{gin}(\mathcal{I}_C)$ is
  obtained from the hyperplane gin in at most two rewritings, by Lemma~\ref{g+i}. Because
  the hyperplane gin is $4$-regular, $\mbox{gin}(\mathcal{I}_C)$, and
  therefore $C$, is $6$-regular. Therefore, $i=0$. It
  follows immediately from \eqref{1.2a} that $C$ is nonproblematic. 

Finally, assume $g \leq 5$. Any $C$ with $g \leq 5$
  that admits an $8$-secant line is nonproblematic, by Lemma~\ref{1a}, and the remaining
curves are $7$-regular, by \cite[Thm.~3.1, p.501]{GLP}. From Figure $7$ it's clear that any
  sequence of at most eight $C$-rewriting rules applied to the hyperplane
  gin that results in a $7$-regular ideal involves at most $3$
  rewritings in degree $6$. Therefore, by Lemma~\ref{4}, $i \leq
  3$. It follows immediately that if $g \geq 4$, then
  $g+i<\min\{2g,8\}$, and therefore such curves are nonproblematic. 

Now say $g \leq 3$. Since curves on hyperquadrics have
codimension 7 in $M_{10}$, by Lemma~\ref{1b}, we may and therefore
shall assume that $\mbox{gin }\mathcal{I}_C$ has no
quadratic generators. Since $\mbox{gin }\mathcal{I}_{\Gamma}$ has {\it two}
quadratic generators, the rewriting sequence that yields
$\mbox{gin}(\mathcal{I}_C)$ from the hyperplane gin necessarily
involves two rewritings in degree $2$. So the generic initial ideal of
$C$ is obtained from the hyperplane gin via a sequence of rewritings
involving at most one rewriting in degree at least six. Therefore, $i
\leq 1$, by Lemma~\ref{4}. It follows that all curves of genus higher
than 1 are nonproblematic. Similarly, if $g=1$ then $g+i<4$ so $C$ is
nonproblematic. Finally, smooth curves are nonproblematic because $i
\leq 1$ and $i \neq 0$ in codimension $4$. 

\begin{Case}\label{C4}
Say $\mbox{gin }\mathcal{I}_{\Gamma}$ admits exactly three
quadratic minimal
generators.
\end{Case}

Then either
\begin{equation*}
\mbox{gin
}\mathcal{I}_{\Gamma}=\mbox{Borel}(x_0x_2)+\mbox{Borel}(x_2^4)+(x_1^3)
\end{equation*}
or
\begin{equation*}
\mbox{gin
}\mathcal{I}_{\Gamma}=\mbox{Borel}(x_0x_1)+(x_0x_2^3,x_1^2,x_1x_2^3,x_2^4).
\end{equation*}
By Lemma~\ref{g_Gamma}, $g_{\Gamma}=9$ in either
case. Therefore, $g+i \leq 9$, by Lemma~\ref{g+i}. 

The proof that the corresponding locus of morphisms is nonproblematic
follows the same lines as in Case~\ref{C3} above, so we merely sketch
it. We may and
shall assume $C$ is $8$-regular. Otherwise by \cite[Thm.~3.1,
  p.501]{GLP}, $C$ has an $8$-secant line; therefore, $C$ belongs to a
proper sublocus of $M_{10}$ of codimension at least $10$, by
Lemma~\ref{1a}. As $g+i<10$, the latter sublocus is nonproblematic.

Next we show that the locus of morphisms
corresponding to $C$ with $g \leq 7$ is nonproblematic. Since we are assuming 
$C$ to be $7$-regular, and since $\mbox{gin}(\mathcal{I}_C$ is obtained from
the hyperplane gin in at most $g_{\Gamma}$ rewritings, it's not hard
to see that at most $5$ rewritings occur in degree $6$, and therefore
$i \leq 5$, by Lemma~\ref{4}. Moreover, our upper bound on $i$
improves to $i \leq 3$ if we also assume $\mbox{gin}(\mathcal{I}_C)$ has no
quadratic minimal generators. If the assumption fails, then the
morphism defining $C$ lies
in a proper sublocus of codimension $7$, by
Lemma~\ref{1b}, but in such a situation $C$ is
nonproblematic. So our second assumption is justified, and we shall
make it. The remainder of the argument when $g \leq 7$ is completely
analogous to the one given in Case~\ref{C3} above.

Finally, say $g=8 \text{ or }9$. The liaison argument given in
Case~\ref{C3} implies that $C$ necessarily lies on a (possibly
singular) cubic scroll. Since rational degree-$10$ curves of genus $8$
on cubic scrolls do not exist, as shown by the analysis of Case~\ref{C3}
above, we may assume $g=9$. Then $C$ is a {\it Castelnuovo curve} in
the sense of \cite[p.27]{Ciliberto} and
\cite[p.527]{GH}; namely, $C$ is of maximal genus among nondegenerate,
irreducible and nondegenerate curves of degree $10$ in $\mathbb{P}^4$. 

Say that $C$ lies on a cubic scroll of type $S_{1,2}$. The
adjunction formula implies that \eqref{adjunction_on_S_{1,2}} holds, with
$g=9$. Solving, we obtain $a=3$ or $a=4$, which correspond to the
classes $3a+7f$ and $4a+6f$, respectively. We next compute
$h^0(\mathcal{O}_{\mathbb{F}_3}(ae+bf)$ for the pairs $(3,7)$ and
$(4,6)$, respectively.

Note that the dimension of the space of cubic scrolls $S_{1,1} \subset
\mathbb{P}^4$ equals $h^0(\mathcal{O}_{\mathcal{F}_1}(e+f))-1$, while the
dimension of the space of curves of class $ae+bf$ on a given scroll equals
$h^0(\mathcal{O}_{S_1}(ae+bf))-1$.

Using the Riemann-Roch formula for surfaces, together with the intersection pairing, we deduce
\begin{equation}
\begin{split}
&\chi(\mathcal{O}_{\mathbb{F}_1}(e+f))=4,\\
&\chi(\mathcal{O}_{\mathbb{F}_1}(3e+7f))=23, \text{ and}\\
&\chi(\mathcal{O}_{\mathbb{F}_1}(4e+6f))=40.\notag
\end{split}
\end{equation}

Moreover, a straightforward cohomological calculation (see
\cite[Section 2]{Cos}) shows that
\begin{equation*}
h^i(\mathcal{O}_{\mathbb{F}_1}(ae+bf))=0, i=1,2,
\end{equation*}
when $(a,b)=(1,1), (3,7), \text{ or }(4,6)$.

If follows that Castelnuovo curves of degree $10$ on cubic scrolls of
type $S_{1,2}$ in
$\mathbb{P}^4$ determine a locally closed sublocus of $M_{10}$ of
dimension at most $44$. On the other hand, each Castelnuovo curve
under our consideration satisfies $i=0$, because each admits a
$6$-regular generic
initial ideal $\mbox{gin}(\mathcal{I}_C)$. Therefore, because
\begin{equation*}
55-44>g,
\end{equation*}
we deduce immediately that the sublocus of $M_{10}$ corresponding to
Castelnuovo curves on scrolls $S_{1,2}$ is nonproblematic.

Finally, say that $C$ lies on a cubic scroll of type $S_{0,3}$. Note
that $C$ cannot pass through the vertex of $S_{0,3}$. To see why, note that adjunction
applied to the proper transform $\stackrel{\sim}{C}$ of $C$ in $\mathbb{F}_3$
yields 
\begin{equation*}
\stackrel{\sim}{g}=7 \text{ and }a-m=3,
\end{equation*}
just as in Case~\ref{C3}, where $a$ is the coefficient of $e$
in the proper transform $[\stackrel{\sim}{C}]$ of $C$ in $\mathbb{F}_{3}$, $m$ is the multiplicity with which $C$
passes through the vertex of the scroll. Because $a-m$ and $a$ are of opposite parity, the 
Riemann-Roch formula implies that either $\chi(\phi_{0,3}^{-1}C)$ or $\chi(\stackrel{\sim}{C})$ is
not an integer, which is absurd.

We conclude that all Castelnuovo curves are nonproblematic. The proof of Theorem~\ref{1.2} is now complete.

\subsection{Curves spanning hyperplanes}
Next, we consider rational curves $C$ whose linear spans are $3$-dimensional
hyperplanes $H \subset \mathbb{P}^4$. By \eqref{1.2b}, in order to
prove Theorem~\ref{1.1} for such curves it suffices to show the
following.
\begin{thm}\label{1.1bisbis}
The reduced irreducible rational curves verifying
\begin{equation*}
g+i \geq 11+ \min(g,5)
\end{equation*}
define a sublocus of $M_{10}$ of codimension greater than $g+i$.
\end{thm}

We now argue much like we did in our analysis of nondegenerate
curves. The analysis here is simpler, though, because the required
bound on $g+i$ is easier to obtain. Therefore, we give the argument,
while omitting the proofs of results that generalize immediately from
the nondegenerate case.

Fix $C$, and let $H_1$ denote the linear span of $C$. Let $H_2$ be
a hyperplane in $H$ that is general with respect to $C$ in
the sense that $\Gamma:=C \cap H_1$ is a collection of ten
points in uniform position in $H$. Choose coordinates for
$\mathbb{P}^4$ in such a
way that $H$ and $H_1$ are defined by $x_4=0$ and $x_3=x_4=0$,
respectively. Just as we did for nondegenerate curves, we define the
{\it hyperplane gin} of $C$ to be the (saturated) generic initial ideal of the
saturation of $\mathcal{I}_{C \cap
  H_1/H}$, and we abusively denote it by
$\mbox{gin}(\mathcal{I}_{\Gamma})$. Just as before, because the hyperplane gin is saturated and
Borel-fixed, it follows that $\mbox{gin }\mathcal{I}_{\Gamma}$ is a monomial ideal in
$x_0$ and $x_1$ having a minimal generating set of the form
\begin{equation*}
(x_0^{k}, x_0^{k-1}x_1^{\lambda_{k-1}}, \dots,
  x_{0}x_1^{\lambda_{1}},x_1^{\lambda_{0}}).
\end{equation*} 

By a result of Ellia and Peskine \cite[Cor.~4.8]{Green}, the {\it invariants} $\lambda_i$ of the above generating set
satisfy
\begin{equation*}
\lambda_i - 1 \geq \lambda_{i+1} \geq \lambda_i - 2
\end{equation*}
for all $i = 0, \dots, k-2$.

Every Borel-fixed ideal with minimal generators that are monomials
$x_0^jx_1^k$ has a unique tree-representation analogous to the tree
representations for the hyperplane gins of nondegenerate curves in
$\mathbb{P}^4$ introduced previously. Moreover, every
  tree may be obtained from an {\it empty tree} $\null$ with a single
  vertex by applying a sequence of rules that we denote as before by
  $\Lambda$-rules (see Table~\ref{Table $3$}). A straightforward inductive
  argument analogous to those already carried out in our analysis of
  nondegenerate hyperplane gins shows that $\sum_i^{k-2} \lambda_i$ is
equal to the number of nonterminal vertices in a minimal generating
tree for the hyperplane gin, which in turn is equal to the degree of
the curve $C$. Whence,
\begin{equation}\label{basiceq}
\sum_i^{k-2} \lambda_i =10.
\end{equation}

On the other hand, by the main result of \cite{Ballico} implies that $I_{\Gamma}$ is 5-regular. It
follows from \eqref{basiceq} and the combinatorial characterization of
Borel-fixity in \cite[Thm.~15.23]{Eisenbud} that the $\mbox{gin}$ of a general hyperplane
section $\Gamma$ of a reduced, irreducible degree-$10$ curve $C$ in $\mathbb{P}^3$ is
either $\mbox{Borel}(x_1^4)$ or $(x_1^5,x_0x_1^3,x_0^2x_1^2,x_0^3)$.
Define $g_{\Gamma}$ for these ideals as before; then a calculation
yields $g_{\Gamma}=11$ and $g_{\Gamma}=12$, respectively. It follows that
\begin{equation}\label{g+i_bound}
g+i \leq 11 \text{ and }
g+i \leq 12,
\end{equation}
respectively.

\begin{table}[section]
\caption{$\Lambda$-rules for curves spanning hyperplanes}\label{Table $3$}
\vspace{-.3cm}
$$\begin{array}{ll}
\hline
\rule{0 pt}{13 pt}1.&x_0^e
\mapsto (x_0^{e+1},x_0^ex_1),\\ 
\rule{0 pt}{13 pt}2.&x_0^ex_1^f \mapsto x_0^ex_1^{f+1},
\text{ and an {\it initial rule}}\\
\rule{0 pt}{13 pt}3.&\emptyset \mapsto
(x_0,x_1,x_2). 
\end{array}$$
\end{table} 

Just as in the case of nondegenerate curves, a minimal generating set
for the generic initial ideal
$\mbox{gin}(\mathcal{I}_{C/H})$ may be represented in a tree
obtainable from the tree of corresponding
hyperplane gin by applying a sequence of rules; these we denote as
before by {\it $C$-rules}. The $C$-rules are given in Table~\ref{Table $4$}.
\begin{table}
\caption{$C$-rules for irreducible curves spanning hyperplanes}\label{Table $4$}
\vspace{-.3cm}
$$\begin{array}{ll}
\hline
\rule{0 pt}{13 pt}1.&x_0^e
\mapsto x_0^e \cdot(x_0,x_1,x_2),\\
\rule{0 pt}{13 pt}2.&x_0^ex_1^f \mapsto x_0^ex_1^f \cdot(x_1,x_2), \text{ and }\\
\rule{0 pt}{13 pt}3.& x_0^ex_1^fx_2^g \mapsto
x_0^ex_1^fx_2^g \cdot x_2
\end{array}$$
\end{table}

The same argument used for nondegenerate curves
shows that for every $C$,
$h^1(\mathcal{I}_{C/H}(5)$ is equal to the number of $C$-rewritings
applied in degrees greater than 6. Finally, the long exact
sequence in cohomology associated to
\begin{equation*}
0 \rightarrow \mathcal{I}_{\mathbb{P}^4}(4) \rightarrow \mathcal{I}_{C/\mathbb{P}^4}(5)
\rightarrow \mathcal{I}_{C/H}(5) \rightarrow 0,
\end{equation*} 
shows that
\begin{equation*}
h^1(\mathcal{I}_{C/H}(5))=h^1(\mathcal{I}_{C/\mathbb{P}^4}(5)).
\end{equation*}

With these preliminaries in hand, the proof of Theorem~\ref{1.1bisbis}
is almost immediate. For, note that applying a single rewriting
rule to either $\mbox{Borel}(x_1^4)$ or
$(x_1^5,x_0x_1^4,x_0^2x_1,x_0^3)$ produces a saturated ideal that is
minimally generated in degrees at most six, and is, therefore,
$6$-regular. Accordingly, our estimates \eqref{g+i_bound} improve to
$g+i \leq 10$ and $g+i \leq
12$, respectively. Therefore, Theorem~\ref{1.1bisbis} is verified in every case except possibly if $g=0$. 

On the other hand, by \cite[Thm.~3.1]{GLP},
$C$ is $8$-regular unless $C$ admits an $8$-secant line. An inspection of possible generic
initial ideals in the second case yields that if $C$ is $8$-regular then
$i \leq 10$. Therefore, Theorem~\ref{1.1bisbis} is verified in every
case away from the locus of genus-$0$ curves that admit $8$-secant
lines. On the other hand, smooth rational curves $C$ of degree $10$ that admit
8-secant lines comprise a sublocus of the Hilbert scheme of rational,
smooth, and irreducible curves in $H$ of codimension at least 5, by
\cite[Lemma (2.4)]{Klei}. In other words, by \eqref{1.2b}, the
sublocus has codimension at least $16$ in $M_{10}$. Theorem~\ref{1.1bisbis} follows
immediately.

\section{Reducible curves}\label{sec-three}
In this section we'll prove the following theorem, extending
\cite[Thm.4.1]{Klei}.

\begin{thm}
On a general quintic threefold in $\mathbb{P}^{4}$, there is no
connected, reduced and reducible curve of degree at most $10$ whose components are rational.
\end{thm}

Suppose, on the contrary, that such a curve $C$ exists. By the results
of \cite{Klei} we may
assume $C$ has two components and is of degree $10$. Consider
one of them. By the result of \cite[Theorem 3.1]{Klei}, either it's a
six-nodal plane quintic or it's smooth. If it's smooth, then, by
\cite[Cor.~2.5(3)]{Klei}, either it's a rational normal curve of
degree $\leq 4$ or it spans
$\mathbb{P}^{4}$. We will prove that there can be
no such $C$.

To this end, we follow Johnsen and Kleiman once more. Let $M_{a}^{'}$ denote the open subscheme of
the Hilbert scheme of $\mathbb{P}^{4}$ parametrizing the smooth irreducible curves of
degree $a$ that are rational normal curves if $a \leq 4$ and that span
$\mathbb{P}^{4}$ if $a \geq 4$. Denote the scheme parametrizing
six-nodal plane quintics in $\mathbb{P}^4$ by $N_{5}$. Let
$R_{a,b,n}$,$S_{5,n}$, and $S_{n}'$ denote the subsets of
$M_{a}^{'} \times M_{b}^{'}$ (resp., $M_5 \times N_5$, $N_{5} \times N_{5}$) of pairs $(A,B)$
such that $A \cap B$ has length $n$. Finally, let $I_{a,b,n}$ (resp., $J_n$,$K_{n}$) denote
the subset of $R_{a,b,n} \times \mathbb{P}^{125}$ (resp., $S_{5,n}
\times \mathbb{P}^{125}$, $S_n' \times \mathbb{P}^{125}$) of triples
$(A,B,F)$ such that $A \subset F$ and $B \subset F$. The $F$ that contain a
plane form a proper closed subset of $\mathbb{P}^{125}$; form its complement, and
the preimages of this complement in the incidence
schemes $I_{a,b,n}, J_n, K_n$. Now replace $J_n$ (resp., $K_n$). Replace $S_{5,n}$ and $S_{n}^{'}$ by the images of those preimages, and
replace $I_{a,b,n}, J_n$ and $K_{n}$ by the new preimages. Then given any pair
$(A,B)$ in $S_{n}^{'}$, there is an $F$ that contains both $A$ and $B$, but
not any plane. It remains to
show that $I_{a,b,n}$ (resp., $J_n$, $K_{n}$)
has dimension at most 124 whenever $a+b = 10$ and $n \geq 1$.

We note that the fibre of $I_{a,b,n}$ (resp., $J_n$, $K_{n}$) over a pair $(A,B)$ is a
projective space of dimension $h^{0}(\mathcal{I}_{C}(5))-1$, where $C$ is the reducible
curve $C=A \cup B$ and $\mathcal{I}_{C}$ is the ideal sheaf of the
corresponding subscheme of $\mathbb{P}^4$. Hence we have

\begin{equation}
\dim I_{a,b,n} \leq \dim R_{a,b,n} + 125 -
\min_{C}\{h^0(\mathcal{O}_C(5))
- h^1(\mathcal{I}_C(5))\} \notag
\end{equation}

Obviously, we have
\begin{equation}
h^{0}(\mathcal{O}_C(5)) \geq \chi(\mathcal{O}_C(5)), \notag
\end{equation}
which implies
\begin{equation}
\begin{split}
&h^0(\mathcal{O}_C(5)) \geq 5(a+b)+2-n, \\
&h^0(\mathcal{O}_C(5)) \geq
20+5a+1-n=21+5a-n, \text{ and}\\
&h^0(\mathcal{O}_C(5)) \geq 20 +
20 - n = 40 - n,\text{ respectively.}\notag
\end{split}
\end{equation}

Our theorem is then a consequence of the following two lemmas.

\begin{Lemma}\label{2.1}
For $a+b=d$ and $n \geq 1$,
\begin{subequations}
\begin{align}
&\dim R_{a,b,n} \leq 5(a+b)+1-n, \label{2.1a}\\
&\dim S_{a,n}
\leq 20+5a-n \label{2.1b}\\ 
\text{and }&\dim S_{n}^{'} \leq
39-n \mbox{.} \label{2.1c} 
\end{align}
\end{subequations}
\end{Lemma}

\begin{Lemma}\label{2.2}
For $a+b=d$ and $n\geq 1$, 
\begin{equation}
h^1(\mathcal{I}_C(5)) = 0. \notag
\end{equation}
\end{Lemma}

To prove the first lemma, begin by letting $(A,B)$ denote an arbitrary
pair in $R_{a,b,n}$. Fix $B$,
and let $A$ vary in the fibre of $R_{a,b,n}$ over $B$. Assume that $a \leq
b$, so in particular $a \leq 4$. 

If $a=1$ or $a=3$, then the lemma holds on the
basis of the arguments in \cite{Klei}. Moreover, if $a=4$ then the
argument of \cite{Klei} carries over except in the case where $B$ is a
sextic meeting $A$ in twelve points along which the sextic intersects
a hyperquadric containing $A$. 

To handle the latter situation, 
recall that the restricted tangent bundles
$T_{\mathbb{P}^4}|_A$ and $T_{\mathbb{P}^4}|_B$ have balanced splittings
    $(5,5,5,5)$ and $(8,8,7,7)$, respectively (\cite[Cor 2.5]{Klei}). Now fix a divisor
$D_1$ of degree 7 along $\mathbb{P}^1$ by and a divisor $D_2$ of degree 7
    along the sextic. The space of degree-6 morphisms $\mathbb{P}^1
      \rightarrow \mathbb{P}^4$ mapping $D_1$ to $D_2$ has dimension
      $h^0(T{\mathbb{P}^4}|_B \otimes \mathcal{I}_{D_2})=30-4\cdot 7$, since
	$h^1(T{\mathbb{P}^4}|_B \otimes \mathcal{I}_{D_2})=0$ (see \cite[p. 45]{Debarre}). As $D_1$ and
      $D_2$ each vary in a 7-dimensional family, it follows that those
	  sextics $B$ intersecting $A$ in seven points cut out a locus
	  of codimension 14 inside $M_{6}$. Since 14 is larger than
	  12, we conclude that $\dim R_{4,6,12}$ meets the
	  required bound.

Similarly, if $a=2$ and $b=8$ then by the argument in \cite{Klei} we may
assume $n \geq 7$ without loss of
generality. Note that if $n \geq 7$, then in fact $n=7$. For, denote
the plane of $A$ by $J$. Note that $B$ spans $\mathbb{P}^{4}$, as $b > 3$. Let $H$ be
the hyperplane spanned by $J$ and a general point of $B$; then $\mbox{l}(H \cap B)
= b \geq n+1$. From $b=8$, we conclude that $n=7$.

To bound the dimension of $R_{2,8,7}$, the space of unions of irreducible
rational conics $A$ and irreducible rational nondegenerate octics $B$
intersecting in projective schemes of length 7, we
proceed much as we did to bound $\dim R_{4,6,12}$. Note it
suffices to show
that rational nondegenerate octics intersecting a conic in at least seven points have
codimension at least 8 inside $M_8$. But just as in the analysis of
$R_{4,6,12}$, the assertion is clear since \cite[Cor 2.5]{Klei}
implies that the
restricted tangent bundle $T_{\mathbb{P}^4}|_B$ has balanced splitting type
  $(10,10,10,10)$.

Next we treat $\dim R_{5,5,n}$, where $n \geq 1$. Using
\cite[Cor 2.5]{Klei}, we note that the 3-regular schemes $A$ and $B$
are each cut out by hyperquadrics and hypercubics.
Therefore,
\begin{equation}
h^0(\mathcal{I}_A(2))=h^0(\mathcal{O}_{\mathbb{P}^4}(2)) - h^0(\mathcal{O}_A(2)) = 15 - (5(2)+1) = 4. \notag
\end{equation}
For degree reasons, no three linearly independent hyperquadrics
containing $A$ may cut out a complete
intersection curve containing $A \cup B$. 

The only remaining possibility is that $A$ and $B$ lie on a
nondegenerate cubic scroll, defined by three
hyperquadrics. But in that case, any ``fourth'' hyperquadric
containing $A$ intersects the scroll properly, in a scheme of degree
at most 6. So we conclude that some hyperquadric $Q$ containing $A$ doesn't
contain $B$. Thus $\mbox{degree }A \cap B \leq \mbox{degree } B \cap Q
=10$. On the other hand, $T_{\mathbb{P}^4}|_A$ has balanced splitting
type $(7,6,6,6)$, from which it follows (by the same argument used
earlier) that for all $1 \leq n \leq
7$, the codimension of curves $A$ intersecting curves $B$ in
subschemes of length $n$ is equal to $2n$. As $\dim R_{5,5,m}
\leq \dim R_{5,5,n}$ whenever $m \geq n$, and $10<14$, it
follows immediately that no curve in $R_{5,5,n}, n\geq 1$ lies on
a general quintic hypersurface.

Now consider reducible unions belonging to $\dim S_{5,n}$, where $n \geq 1$. Let $A$ be a
nondegenerate, rational, smooth quintic curve, and let $B$ be a
six-nodal plane quintic with linear span $J$. The intersection $A \cap
J$ is proper, whence of degree at most 5; it follows that $\deg(A \cap
B) \leq 5$. Since $T_{\mathbb{P}^4}|_A$ has splitting
type $(7,6,6,6)$, the codimension of curves $A$ intersecting curves $B$ in
subschemes of length $n$ is equal to $2n$, for all $n \leq 7$. It
follows that no curve in $S_{5,n}, n \geq 1$ lies on a general quintic hypersurface.

Finally, consider a pair $(A,B)$ in $S_{n}^{'}$. Let $J$ denote
the plane of $A$, and $K$ the plane of $B$. We may assume $J \neq K$
without loss of generality, since the general quintic threefold $F$
intersects $J$ properly in a quintic curve, which shows that if $J=K$, then
\begin{equation*}
A \cup B \subset F
\end{equation*}
is impossible. So $A \cup B$ spans at least a 3-space. 

Moreover, it's clear that $A \cap B \subset J \cap B$, so $A
\cap B$ has degree at most 5. On the other hand, the restricted
tangent bundle $T_{\mathbb{P}^2}|_A$ has splitting type $(a_1,a_2)$,
where $a_1+a_2=15$. Assume $a_1 \geq a_2$. As usual, our goal is to
bound the codimension of those six-nodal plane quintics that meet
other six-nodal plane quintics in $1 \leq n \leq 5$ points using
ampleness properties of the restricted tangent bundle of $A$. In
particular, we are done provided $a_2 \geq 2$, which is certainly the
case. For (cf. \cite{GLP}), there is an exact sequence
\begin{equation}
0 \rightarrow \bigoplus_{i=1}^2 \mathcal{O}_{\mathbb{P}^1}(-a_i+5) \rightarrow
H^0(\mathcal{O}_{\mathbb{P}^1}(5)) \otimes \mathcal{O}_{\mathbb{P}^1} \rightarrow
\mathcal{O}_{\mathbb{P}^1}(5) \rightarrow 0\notag
\end{equation}
which implies $a_i-5 \geq 0$.
The proof of the first lemma is now complete.  

We now proceed to the proof of the second lemma, i.e. that all curves
$C$ in rational components $A$ and $B$ satisfying our hypotheses are
6-regular. First say $A$ and $B$ are smooth. It's a well-known fact
(see, for example, \cite[Thm 2.1]{Gi}) that
\begin{equation*}
\mbox{reg }A \cup B \leq \mbox{reg }A+ \mbox{reg }B,
\end{equation*}
so we need only establish that $A$ and $B$ are 3-regular. But this
follows, e.g., from the result of \cite[Prop. 2.2]{Klei2}. (See also
the discussion following the proof of \cite[Cor.~2.5]{Klei}; the key
point is to observe that the components $C_i$ of $C$ are necessarily of
{\it maximal rank} in every degree for the canonical morphisms
$H^0(\mathcal{O}_{\mathbb{P}^4}(k)) \rightarrow H^0(\mathcal{O}_{C_i}(k))$.) 

In the only remaining case, $A$ is a smooth rational quintic spanning $\mathbb{P}^4$
and $B$ is a six-nodal plane quintic. In that situation, $C=A \cup B$ is
$4$-regular, by the ``Horace lemma'' \cite[Lemma 4.5]{Klei}.

Therefore, we may assume $A$ and
$B$ are plane quintics. Once more, we treat separately the cases
where $\mbox{rank }C=3$ and $\mbox{rank }C=4$.

If $G$, the linear span of $A \cup B$, is a 3-space, then $G \cap H$ is a plane containing $K$,
and is therefore equal to $K$. Hence 
\begin{equation}
C \cap H = (C \cap G) \cap H = C \cap (G \cap H) = C \cap K. \notag
\end{equation}
On the other hand, we have
\begin{equation}
(A \cup B) \cap K = (A \cap K) \cup B. \notag
\end{equation}

Indeed, the latter equality of schemes is equivalent to the following equality
of ideal sheaves:
\begin{equation}
(\mathcal{I}_{A/\mathbb{P}^4} \cap \mathcal{I}_{B/\mathbb{P}^4}) + \mathcal{I}_{K/\mathbb{P}^4} =
 (\mathcal{I}_{A/\mathbb{P}^4}+\mathcal{I}_{K/\mathbb{P}^4}) \cap \mathcal{I}_{B/\mathbb{P}^4}. \notag
\end{equation}
Any element $l$ belonging to the left side is of the form $l=a+k=b+k$ where
$a \in \mathcal{I}_{A/\mathbb{P}^4}$, $b \in \mathcal{I}_{B/\mathbb{P}^4}$, and $k \in \mathcal{I}_{K/\mathbb{P}^4}$. The
inclusion $l \in \mathcal{I}_{A/\mathbb{P}^4}+\mathcal{I}_{K/\mathbb{P}^4}$ follows immediately, and we also
have $l \in \mathcal{I}_{B/\mathbb{P}^4}$ because $\mathcal{I}_{K/\mathbb{P}^4} \subset \mathcal{I}_{B/\mathbb{P}^4}$ from the
inclusion $B \subset K$.  In the opposite direction, given $r = a+k = b$,
we see that $a = b-k \in \mathcal{I}_{B/\mathbb{P}^4}$, again from $B \subset
K$.  Since a general quintic threefold $F$ contains no plane, we may
assume $A \cup B \subset F$ but $K$ lies outside $F$, so $F \cap K=B$
by Bezout's theorem. Then 
$A \cap K \subset B$, and it follows immediately that
$C \cap H = C \cap K = B$.

In what follows, we let $D$ denote $C \cap H$.
To compute $h^1(\mathcal{I}_{C/\mathbb{P}^4}(5))$, we use the exact sequence:
\begin{equation}
0 \rightarrow \mathcal{I}_{A/\mathbb{P}^4}(-1) \rightarrow \mathcal{I}_{C/\mathbb{P}^4}
\rightarrow \mathcal{I}_{D/H} \rightarrow 0. \notag
\end{equation}

To bound the cohomology of the middle term, we bound cohomology on the
right and left.
Note there is an exact sequence
\begin{equation}
0 \rightarrow \mathcal{I}_{K/H} \rightarrow \mathcal{I}_{D/H} \rightarrow \mathcal{I}_{D/K}
\rightarrow 0 \notag
\end{equation}
with $\mathcal{I}_{K/H} = \mathcal{O}_H(-1) = \mathcal{O}_{\mathbb{P}^3}(-1)$ and $\mathcal{I}_{D/K} = \mathcal{O}_K(-5) =
\mathcal{O}_{\mathbb{P}^2}(-5)$. It follows immediately from Serre's theorem that
$h^1(\mathcal{I}_{K/H}(m)) = h^1(\mathcal{I}_{D/K}(m)) = 0$ for all integers $m$.
Similarly, there is an exact sequence
\begin{equation}
0 \rightarrow \mathcal{I}_{J/\mathbb{P}^4} \rightarrow \mathcal{I}_{A/\mathbb{P}^4}
\rightarrow \mathcal{I}_{A/J} \rightarrow 0 \notag
\end{equation}
with $\mathcal{I}_{A/J} = \mathcal{O}_{\mathbb{P}^2}(-5)$ and $\mathcal{I}_{J/\mathbb{P}^4} =
\mathcal{O}_{\mathbb{P}^4}(-1)^{2}$, and hence $h^1(\mathcal{I}_{A/\mathbb{P}^4}(m)) = 0$ for all integers $m$.

Finally, say $C= A \cup B$ is nondegenerate. After having made an appropriate change of
coordinates, we may assume that $J$ and $K$ intersect in the point $P =
(0,0,1,0,0)$ and that the homogeneous ideals describing the embeddings
of $A$ and $B$ inside $\mathbb{P}^4$ are given in coordinates by
\begin{equation}
I_A = (x_3, x_4, f(x_0,x_1,x_2)) \mbox{ and } I_B = (x_0, x_1,
g(x_2,x_3,x_4)) \notag
\end{equation}
for some trivariate homogeneous quintic polynomials $f$ and $g$. Since $A$
and $B$ pass through $P = (0,0,1,0,0)$, $f$ and $g$, which vanish at $P$,
do not contain $x_2^5$ in their expansions.

Note that
\begin{equation}
(x_3, x_4) \cap I_B + f(x_0,x_1,x_2) \cap I_B \subset I_A \cap I_B \notag.
\end{equation}
On the other hand, given any element $e \in I_A \cap I_B$, viewed as a
combination of $x_3$, $x_4$, and $f(x_0,x_1,x_2)$ with polynomial coefficients, any terms of $e$ divisible by $f(x_0,x_1,x_2)$ automatically belong to $I_B$, since $f$ does not contain $x_2^5$ in its expansion. It follows immediately that $e=e_1+e_2$ with $e_2 \in f(x_0,x_1,x_2) \cap I_B$ and $e_1 \in (x_3, x_4) \cap I_B$, and therefore that
\begin{equation}
(x_3, x_4) \cap I_B + (f(x_0,x_1,x_2)) \cap I_B = (x_3, x_4) \cap I_B
  + (f(x_0,x_1,x_2)) = I_A \cap I_B.  \notag
\end{equation}

Continuing in this vein, we deduce that the homogeneous ideal
of $C \subset \mathbb{P}^4$ is 
\begin{equation}
\begin{split}
&I_C = I_A \cap I_B = (x_3, x_4) \cap (x_0, x_1) + (g(x_2,x_3,x_4)) + (f(x_0,x_1,x_2)) \\
&= (x_1x_4,x_0x_4,x_1x_3,x_0x_3,f(x_0,x_1,x_2),g(x_2,x_3,x_4)). \notag 
\end{split}
\end{equation}
We will now show that 
\begin{equation}
(x_1x_4,x_0x_4,x_1x_3,x_0x_3,\mbox{lt}(f),
\mbox{lt}
(g)) = \mbox{in}(I_C), \notag
\end{equation}
where $\mbox{lt}(F)$ denotes the leading term of a homogeneous
polynomial $F$ with respect to the revlex order, and
$\mbox{in}(I)$ denotes the ideal of leading terms of the homogeneous
ideal $I$.

We will deduce the latter inclusion as a consequence of \cite[Prop 4.3]{Green},
which establishes that for any two homogeneous ideals $I$ and $J$, if
\begin{equation}
\mbox{in}(I) \cap \mbox{in}(J) \subset \mbox{in}(I \cap J) \notag
\end{equation}
then $\mbox{in}(I+J)=\mbox{in}(I)+\mbox{in}(J)$.

Now let 
\begin{equation}
f=\sum_{i+j+k=5}a_{ijk} x_0^ix_1^jx_2^k; \notag
\end{equation} 
since $(i,j)=(0,0)$
is disallowed, we have 
\begin{equation}
x_3f \in (x_1x_4,x_0x_4,x_1x_3,x_0x_3)
  \cap (f). \notag
\end{equation}

So 
\begin{equation*}
\mbox{lt}(x_3f)=x_3\mbox{lt}(f) \in \mbox{in}((x_1x_4,x_0x_4,x_1x_3,x_0x_3)
  \cap (f))
\end{equation*} and similarly 
\begin{equation*}
\mbox{lt}(x_4f) \in \mbox{in}((x_1x_4,x_0x_4,x_1x_3,x_0x_3)
  \cap (f)).
\end{equation*}
On the other hand, if 
\begin{equation*}
\mbox{lt}(f)=a_{ijk}x_0^ix_1^jx_2^k,
\end{equation*}
then $x_3a_{ijk}x_0^ix_1^jx_2^k$
and
  $x_4a_{ijk}x_0^ix_1^jx_2^k$ generate the intersection 
\begin{equation}
\mbox{in}(x_1x_4,x_0x_4,x_1x_3,x_0x_3) \cap \mbox{in}(f). \notag
\end{equation}
 (This is clear, since any monomial $e$ belonging to the latter intersection
  is divisible by one of the monomials $x_1x_4,x_0x_4,x_1x_3,x_0x_3$,
  hence by $x_3$ or $x_4$, and also $\mbox{lt }f$. Thus $e$ is divisible
  by $x_3
  \mbox{lt }f$ or $x_4\mbox{lt }f$.) So in fact  
\begin{equation}
\mbox{in}(x_1x_4,x_0x_4,x_1x_3,x_0x_3) \cap \mbox{in}(f) \subset \mbox{in}((x_1x_4,x_0x_4,x_1x_3,x_0x_3)
  \cap (f)). \notag
\end{equation}
Therefore,
\begin{equation}
\mbox{in}(x_1x_4,x_0x_4,x_1x_3,x_0x_3,f)=(x_1x_4,x_0x_4,x_1x_3,x_0x_3,\mbox{lt}(f)).\notag
\end{equation}

Essentially the same argument shows that
\begin{equation}
\mbox{in}(x_1x_4,x_0x_4,x_1x_3,x_0x_3,f)\cap \mbox{in}(g) \subset \mbox{in}((x_1x_4,x_0x_4,x_1x_3,x_0x_3,f)
  \cap (g)). \notag
\end{equation}
We conclude that
\begin{equation}
\mbox{in}(x_1x_4,x_0x_4,x_1x_3,x_0x_3,f,g)=(x_1x_4,x_0x_4,x_1x_3,x_0x_3,\mbox{lt}(f),\mbox{lt}(g)). \notag
\end{equation}

Since there is a flat degeneration taking $\mathcal{I}_C$ to 
$\mbox{in }\mathcal{I}_C$, and the latter ideal is generated in degrees  at most
5,
we conclude by the basic regularity result of Bayer--Stillman
\cite[Theorem 2.27]{Green} that $\mathcal{I}_C$ is $6$-regular.

\end{document}